\documentclass[a4paper,12pt]{article}
\usepackage[T2A]{fontenc}             
\usepackage[utf8]{inputenc}           
\usepackage[russian,english]{babel}          
\usepackage{amsmath}                  
\usepackage{amsfonts}
\usepackage{amssymb}
\usepackage{exscale}                 
\usepackage[top=2cm, left=2cm, right=2cm, bottom=2cm]{geometry} 
\numberwithin{equation}{section}
\usepackage[pdftex,unicode]{hyperref}  
\hypersetup{colorlinks, citecolor=blue, filecolor=blue, linkcolor=blue, urlcolor=blue}                    
\newtheorem{Lemma}{Лемма}[section]
\newtheorem{Theorm}{Теорема}[section]

\begin{document}
\title{Expansion into a many-dimensional rational series for scalar power  functions of vector arguments.}
\author{Robert Akhmetyanov F.\footnote{robertu@mail.ru} , Elena Shikhovtseva S.  \\
	 \itshape Institute of Molecule and Crystal Physics URC RAS \\  
	 \itshape (Prospekt Oktyabrya 71, Ufa, Russia, 450054) }
\date{}	
\maketitle         

\begin{abstract}
For a function of a type $ \left| \mathbf{r}_1{+}\ldots {+}\mathbf{r}_{_N} \right|^{-\nu} \in \mathbb{R} $ from the many-dimensional vectors $ \mathbf{r}_s $ in Euclidean space, the successive algebraic approach is the derivation of the expansion in the form  $ {\sim}\sum\limits_{s}\frac{r_1^s}{r_{_{N}}^s}\ldots \frac{r_{_{N-1}}^s}{r_{_{N}}^s}, \,\, (r_k{<}r_{_{N}}) $, and also for certain orthogonal functions $ H_{\lambda_s}(\mathbf{r}_s) $ as $ {\sim}\sum\limits_{\lambda_k} H_{\lambda_1}(\mathbf{r}_1)\ldots H_{\lambda_{_{N}}}(\mathbf{r}_{_{N}}) $. The coefficient angular functions are found and determined in both cases.
\end{abstract}
\begin{quote}
\end{quote}

\newpage
\selectlanguage{russian}
\begin{center}
{\Large \textbf{Разложение в многомерный рациональный ряд для скалярных степенных функций векторных  аргументов.}} \\
\bigskip
Ахметьянов Р.Ф.,  Шиховцева Е.С. \\
\medskip
\textit{ Институт физики молекул и кристаллов УНЦ РАН, \\ 
Россия, 450054, г.Уфа, Пр.Октября, 71 \\
E-mail: robertu@mail.ru}
\end{center}

\begin{quote}
 \textbf{Аннотация:} Для функции вида $ \left| \mathbf{r}_1{+}\ldots {+}\mathbf{r}_{_N} \right|^{-\nu} {\in} \mathbb{R} $ от многомерных векторов $ \mathbf{r}_s $ в евклидовом пространстве, последовательным алгебраическим подходом представлен вывод разложения в виде $ {\sim}\sum\limits_{s}\frac{r_1^s}{r_{_{N}}^s}\ldots \frac{r_{_{N-1}}^s}{r_{_{N}}^s}, \,\, (r_k{<}r_{_{N}}) $, а также по определенным ортогональным функциям $ H_{\lambda_s}(\mathbf{r}_s) $ в виде  $ {\sim}\sum\limits_{\lambda_k} H_{\lambda_1}(\mathbf{r}_1)\ldots H_{\lambda_{_{N}}}(\mathbf{r}_{_{N}}) $. Найдены и определены коэффициентные угловые функции в обоих случаях. \\

 \textbf{Ключевые слова:} обобщенная гипергеометрическая функция, гиперсферическая функция, функция Гегенбауэра.
\end{quote}

\section{Введение и обозначения} \label{I}
Применение разложения скалярной функции от трехмерных векторов $ \mathbf{r}_s{=}r_s \boldsymbol{\zeta}_s, \, |\mathbf{r}_s|{=}r_s $  по сферическим функциям $ Y_{l,m}(\boldsymbol{\zeta}_s) $

\begin{equation}\label{1_a0}
\begin{gathered}
\frac{1}{|\mathbf{r}_1-\mathbf{r}_2|}
=4\pi\sum_{l=0}^{\infty}\sum_{m=-l}^{l}
\frac{1}{2l+1}\frac{r_{<}^l}{r_{>}^{l+1}}
Y_{l,m}(\boldsymbol{\zeta}_1)Y_{l,m}^{\ast}(\boldsymbol{\zeta}_2) 
\\
r_{<}=\min(r_1,r_2), \quad r_{>}=\max(r_1,r_2)
\end{gathered}
\end{equation}
широко используется в физических и математических задач, обладающих сферической симметрии.
Однако возможно особый интерес в физических и математических приложениях и задачах представляет не разделение по отдельности угловым  $ \boldsymbol{\zeta}_s $ и пространственным $ r_s $  переменным как в \eqref{1_a0}, а разделение по векторам $ \mathbf{r}_s $. Как было показано в \cite{AkhR_UNC} такое разделение существует для трехмерных двух векторов, (для всех действительных $ \nu<3 $)
\begin{multline}
\label{1_a015}
\frac{1}{\left| \mathbf{r}_1 - \mathbf{r}_2 \right|^{\nu}}=
\frac{\pi^{\tfrac{3}{2}}\Gamma\left(\dfrac{3-\nu}{2}\right)}{\Gamma\left( \dfrac{\nu}{2}\right)} 
\left(\left(r_1^2+1\right) \left(r_2^2+1\right) \right)^{\tfrac{3-\nu}{2}} 
\times 
\\
\times 
\sum_{n=0}^{\infty} \sum_{l=0}^{\infty} \sum_{m=-l}^{l}
\frac{\Gamma\left(l+n+\dfrac{\nu}{2} \right)}{\Gamma\left(l+n+3-\dfrac{\nu}{2}\right)} 
H_{n,l,m}(\mathbf{r}_1) 
H_{n,l,m}(\mathbf{r}_2)^{\ast}
\end{multline}
где
\begin{equation} \label{1_eH3}
H_{n,l,m}(\mathbf{r})=\eta_{n,l}(r)Y_{l,m}(\boldsymbol{\zeta})
\end{equation}
\begin{multline} \label{1_e_eta3}
\eta_{n,l}(r)=
\frac{2}{\Gamma\left(l+\dfrac{3}{2}\right)}
\sqrt{\frac{\left(n+l+1\right)\Gamma\left(n+2l+2\right)}{\Gamma\left(n+1\right)}} 
\times 
\\
\times 
\frac{r^l}{\left(r^2+1\right)^{l+\tfrac{3}{2}}}
\,{}_2F_1 \left[ \left. 
\begin{matrix}
{  -n  \quad n+2l+2  } \\
{  l+\dfrac{3}{2}  }
\end{matrix} \right| \frac{1}{r^2+1}\right] 
\end{multline}

Здесь мы рассмотрим аналогичное разложение \eqref{1_a015}, но для произвольных $ M $ ---мерных векторов, и с любым $ N $  их количеством в евклидовом пространстве.

Везде и в дальнейшем  $ C_{\mu}^{\alpha}(z) $ ---функция Гегенбауэра.
Под символами $ (a)_k $  будет означаться символ Похгаммера
\[ (a)_k=a(a+1)\ldots (a+k-1),\,\,[k=1,2,3,\ldots], \quad (a)_0=1 \]
или через Гамму функцию $ \Gamma $   
\[ (a)_k=\frac{\Gamma( a+k )}{\Gamma( a )} \]
Обобщенная гипергеометрическая функция есть как
\[ \,{}_{p}F_{q} \left[ 
\left.
\begin{gathered}
{ (a_p) } \\
{ (b_q) }
\end{gathered}
\right|  z  \right]=\sum_{k=0}^{\infty}\frac{\prod\limits_{p}(a_p)_k}{\prod\limits_{q}(b_q)_k} \frac{z^k}{k!} \]
\[ (a_p)=a_1,a_2,\ldots ,a_p, \,\, \prod_{p}(a_p)_k=\prod_{j=1}^{p}(a_j)_k=\prod_{j=1}^{p}\frac{\Gamma( a_j+k )}{\Gamma( a_j )} \]

Многомерные единичные вектора в $ M $ -мерном евклидовом пространстве мы будем обозначать как $ \boldsymbol{\zeta}_s $ , где  $ s $-означает принадлежность (номер) вектора. При этом скалярное произведение, в силу единичности векторов есть как косинус угла между ними $ \cos \omega_{\alpha \beta} =(\boldsymbol{\zeta}_{\alpha}\boldsymbol{\zeta}_{\beta}) $. Мы не будем конкретизировать выбор полярных координат и их компонентов, всего их может существовать $ \dfrac{(2M-2)!}{(M-1)!\,M!} $  эквивалентных представлении. Соответственно столько же эквивалентных представлении существует и для гиперсферических функции $ Y_{l,\mathbf{m}_k^{(s)}}(\boldsymbol{\zeta}_s) $,  $ \mathbf{m}_k^{(s)}=\{ l_s=m_0^{(s)},m_1^{(s)},\ldots,m_{_{M-2}}^{(s)} \}  $  ( здесь $ s $  также означает принадлежность индексов к соответствующему единичному вектору). Теория гиперсферических функции хорошо изложена во многих литературах (к примеру \cite[Гл.11]{Beitman_2}). При этом будем полагать такой выбор этих функции, при котором
\begin{gather} \label{int_Y}
\int\!\! d\Omega_{\boldsymbol{\zeta}}Y_{l_1,\mathbf{m}_k^{(1)}}(\boldsymbol{\zeta}) Y_{l_2,\mathbf{m}_k^{(2)}}(\boldsymbol{\zeta})^{\ast}=\delta_{l_1,l_2}\delta_{\mathbf{m}_k^{(1)},\mathbf{m}_k^{(2)}} \\ \delta_{\mathbf{m}_k^{(1)},\mathbf{m}_k^{(2)}}=\delta_{m_1^{(1)},m_2^{(2)}}\delta_{m_2^{(1)},m_2^{(2)}}\ldots \delta_{m_{_{M-2}}^{(1)},m_{_{M-2}}^{(2)}} 
\end{gather}
где $ \delta_{p_1,p_2} $---символ Кронекера. Во всех разделах этой работы интеграл по поверхности единичной гиперсферы  $ \int\!\! d\Omega_{\boldsymbol{\zeta}}f(\boldsymbol{\zeta}) $ подразумевается, что интегрирование берется по всему $ M-1 $ мерному пространству. Для произвольной системы гиперсферических координат, элемент объема и его площадь $ M{-}1 $ --мерной сферы представляются соотношениями соответственно как \cite[Гл.11]{Beitman_2}
\[ d\mathbf{r}=dV_{_{M}}=r^{M-1}dr\,d\Omega_{\boldsymbol{\zeta}},\,\, S_{_{M}}=\frac{2\pi^{\tfrac{M}{2}}}{\Gamma\left( \dfrac{M}{2} \right)} \]

Отметим, что при $ M=3 $ для \eqref{1_eH3}  выполняется условие ортогональности
\begin{equation*}
\int\!\! d\mathbf{r} \, H_{n_1,l_1,m_1}(\mathbf{r})H_{n_2,l_2,m_2}(\mathbf{r})^{\ast}=
\delta_{n_1,n_2}\delta_{l_1,l_2}\delta_{m_1,m_2}
\end{equation*}

\section{Представление вида  \eqref{1_a0} для  $ N>2 $ многомерных векторов} \label{sec_1_2}

\begin{Theorm} \label{Th_1_2}
	Для произвольных $ M $ ---мерных векторов $ \mathbf{r}_{s}=r_s \cdot \boldsymbol{\zeta}_{s} \quad s=1,2,\ldots,N $, где $ \boldsymbol{\zeta}_{s} $ единичные $ M $ ---мерные вектора в евклидовом пространстве, а также для любых  $ \nu \in \mathbb{R} $ для однородных функции  $ \left| \mathbf{r}_1{+}\ldots {+}\mathbf{r}_{_N} \right|^{-\nu} \in \mathbb{R} $ справедливо выражение   
	\begin{multline} \label{1_b12}
	\frac{1}{\left|\mathbf{r}_1+\ldots+\mathbf{r}_{_N}\right|^{\nu}}= 
	\frac{1}{(r_{_N})^{\nu}} \sum_{l_1,\ldots,l_{_N}=0}^{\infty} (-1)^{l_{_N}}  V_{l_1,\ldots, l_{_N}} \left( \boldsymbol{\zeta}_{1},\ldots ,\boldsymbol{\zeta}_{N} \right)
	\times \\ \times 	
	\sum_{\mu_1,\ldots,\mu_{_{N-1}}=0}^{\infty}
	\frac{\left(\dfrac{\nu}{2}\right)_{\mu_1+\ldots +\mu_{_{N-1}}+l} \left(\dfrac{\nu{-}M{+}2}{2}\right)_{\mu_1+\ldots +\mu_{_{N-1}}+l-l_{_N}} }
	{\displaystyle\prod\limits_{p=1}^{N-1} \left[\left(\dfrac{M}{2}\right)_{l_p+\mu_p} \mu_p ! \right]} 
	\prod\limits_{p=1}^{N-1} \left(\frac{r_p}{r_{_N}}\right)^{l_p+2\mu_p}
	\end{multline} 
	где
	\[ l=\frac12 \sum_{p=1}^{N} l_p , \quad \max (r_1,r_2,\ldots ,r_{_N})=r_{_N} \] 
	а коэффициентная угловая функция выражается через  гиперсферические \\  функции $ Y_{l,\mathbf{m}_{\mathbf{k}} } \left(\boldsymbol{\zeta} \right) $ или через функции Гегенбауэра от скалярных произведении \\ единичных векторов $ C_{l_i}^{\tfrac{M}{2}-1}\left((\boldsymbol{\zeta}_{i} \boldsymbol{\zeta})\right)  $   как
	\begin{multline} \label{1_b11}
	V_{l_1,\ldots, l_{_N}} \left( \boldsymbol{\zeta}_{1},\ldots ,\boldsymbol{\zeta}_{N} \right)=
    \int \frac{d\Omega_{\boldsymbol{\zeta}}}{S_{_M}} \prod\limits_{i=1}^{N} \frac{l_i+\dfrac{M}{2}-1}{\dfrac{M}{2}-1} 
	C_{l_i}^{\tfrac{M}{2}-1}\left((\boldsymbol{\zeta}_{i} \boldsymbol{\zeta})\right)= 
	\\
	=\left(S_{_M}\right)^{N-1}\sum_{\mathbf{m}_{\mathbf{k}}^{(1)},\ldots,\mathbf{m}_{\mathbf{k}}^{(N)} } \,\,\, 
	\prod_{p=1}^{N} Y_{l_p,\mathbf{m}_{\mathbf{k}}^{(p)} } \left(\boldsymbol{\zeta}_{p} \right) 
	\int d\Omega_{\boldsymbol{\zeta}} \prod_{p=1}^{N} \left(Y_{l_p,\mathbf{m}_{\mathbf{k}}^{(p)} } \left(\boldsymbol{\zeta} \right)\right)^{\ast} 
	\end{multline} 
	
\end{Theorm}

Выбор формы записи коэффициентой угловой функции обусловлен тем, что при допустим $ l_n=0 $  уменьшается порядок  $ N $. Например, при $ N=5 $ и $ l_1=0,\,l_4=0 $ получим
\[  V_{0,l_2,l_3,0,l_5} \left( \boldsymbol{\zeta}_{1},\boldsymbol{\zeta}_{2},\boldsymbol{\zeta}_{3},\boldsymbol{\zeta}_{4},\boldsymbol{\zeta}_{5} \right) =
V_{l_2,l_3,l_5} \left( \boldsymbol{\zeta}_{2},\boldsymbol{\zeta}_{3},\boldsymbol{\zeta}_{5} \right) \]
то есть мы получили порядок $ N=3 $.
Также следует отметить, что \eqref{1_b11} отлично от нуля когда $ \sum\limits_{p=1}^{N}l_p $ являются четным числом, а следовательно $ l $-- целое число. При $ \nu=-2q, \,\, q=1,2,3,\ldots $ суммы в \eqref{1_b12} конечны, так как при $ \left( -q \right)_{\mu_1+\ldots+\mu_{_{N-1}}+l}\neq 0 $ (даже если $ M $--четно) приводит к ограничениям по $ \mu_p $ и $ l_p $ как $ \mu_1+\ldots+\mu_{_{N-1}}+l \leqslant q $ или $ l_{_N}+\sum\limits_{p=1}^{N-1} (l_p+2\mu_p) \leqslant 2q $.
\\
Отметим, что значение $  \max (r_1,r_2,\ldots ,r_{_N})=r_{_N} $   условно в наших обозначениях, так как из условия коммутативности левого выражения в \eqref{1_b12} мы всегда можем выбрать  такой индекс $ N $, путем перестановки или переобозначении. Для краткости обозначении сумм здесь и везде будем полагать что
\[
\sum_{\lambda_1,\ldots,\lambda_k=0}^{\infty}=\sum_{\lambda_1=0}^{\infty}\ldots \sum_{\lambda_k=0}^{\infty} 
\]
\\ 
Выражение \eqref{1_b12} по сути является аналогичным выражением \eqref{1_a0} при $ b=\dfrac{3}{2} $, $ N=2 $ и $ \nu=1 $

Для получения выражения \eqref{1_b12} и \eqref{1_b11} будем исходить из производящей функции для полиномов Гегенбауэра, в котором в нашем случае запишем как
\begin{multline}
\label{1_b2}
\frac{1}{\left|\mathbf{R}_{_{N-1}}+\mathbf{r}_{_N}\right|^{2\rho}}=
\frac{1}
{(r_{_N})^{2\rho}\left(1+2\dfrac{R_{_{N-1}}}{r_{_N}} (\boldsymbol{\zeta}_{R} \boldsymbol{\zeta}_{N})+\left(\dfrac{R_{_{N-1}}}{r_{_N}}\right)^2 \right)^{\rho} }=
\\
=\frac{1}{(r_{_N})^{2\rho}} \sum_{l=0}^{\infty}\left(-\frac{R_{_{N-1}}}{r_{_N}}\right)^l C_l^{\rho}\left((\boldsymbol{\zeta}_{R} \boldsymbol{\zeta}_{N})\right)=
\\ 
=\frac{1}{(r_{_N})^{2\rho}} \sum_{\sigma=0}^{1} \sum_{l=0}^{\infty}\left(\frac{R_{_{N-1}}}{r_{_N}}\right)^{2l+\sigma} (-1)^{\sigma} C_{2l+\sigma}^{\rho}\left((\boldsymbol{\zeta}_{R} \boldsymbol{\zeta}_{N})\right)
\end{multline}
где $  \mathbf{R}_{_{N-1}}{=}R_{_{N-1}}\boldsymbol{\zeta}_{R}, \,\, \rho \in \mathbb{R}  $. Отметим, что в \eqref{1_b2} $ \dfrac{R_{_{N-1}}}{r_{_N}}<1 $, а в последнем выражении мы разложили сумму по четным ($ \sigma=0 $) и нечетным  ($ \sigma=1 $) суммам. 

Для начала покажем, что для полиномов Гегенбауэра справедлива лемма
\begin{Lemma} \label{Lm_1_2}
	Для любых единичных $ M $--мерных векторов векторов $ \boldsymbol{\zeta}_{1} $ и $ \boldsymbol{\zeta}_{2} $ справедливо тождество	
	\begin{multline}
	\label{1_b3}
	C_{2l+\sigma}^{\rho}\left((\boldsymbol{\zeta}_{1} \boldsymbol{\zeta}_{2})\right)=
	\\
	=
	\sum_{m=0}^{l}\frac{(-1)^{m+l}\Gamma(\rho+2l+\sigma+1)\Gamma(\rho+l+m+\sigma)\Gamma\left(\dfrac{M}{2}+l+m+\sigma\right) }{(l-m)! \,\, \Gamma(\rho)\Gamma(\rho+l+m+\sigma+1)\Gamma\left(\dfrac{M}{2}\right) } 
	\times 
	\\
	\times \int\frac{d\Omega_{\boldsymbol{\zeta}}}{S_{_M}} \frac{\left(2(\boldsymbol{\zeta}_{1} \boldsymbol{\zeta})\right)^{2l+\sigma}}{(2l+\sigma)!} \frac{\left(2(\boldsymbol{\zeta}_{2} \boldsymbol{\zeta})\right)^{2m+\sigma}}{(2m+\sigma)!}
	\end{multline}
	где $ \sigma= 0  $ или $ \sigma=1 $. 
\end{Lemma}
Действительно, представим функцию Гегенбауэра в виде  ~\cite[Гл.3, п.3.15.1]{Beitman_1} \\ (где $ \sigma=0,1 \quad l=0,1,2,\ldots $)
\begin{multline}
\label{1_b4}
C_{2l+\sigma}^{\rho}(z)=
\frac{(-1)^l\Gamma(\rho+l+\sigma) }{\Gamma(l+1)\Gamma(\rho)}(2z)^{\sigma}
{}_2F_1 \left[ \left. \begin{gathered}
{-l \quad \rho+l+\sigma} \\
{\sigma+\dfrac12}
\end{gathered} \right| z^2\right] = \\
=\sum_{k=0}^{l} \frac{(-1)^{l+k}\Gamma(\rho+l+\sigma+k)}{\Gamma(l-k+1)\Gamma(\rho)} \frac{(2z)^{2k+\sigma}}{(2k+\sigma)!}
\end{multline}
из выражения \cite[Гл.5 п.5.3.1(9)]{Prudnikov_3} для конечных сумм от обобщенных гипергеометрических функции (где $ \binom{l}{m}=\frac{l!}{m! \,(l-m)!} $)
\begin{equation*}
\sum_{m=0}^{l} \frac{(-1)^m \, l! \, (\alpha)_m}{m! \, (l-m)! \, (\beta)_m} 
{}_{p+1}F_q \left[ \left. \begin{gathered}
{-m \quad (a_p)} \\
(b_q)
\end{gathered} \right| x \right]=
\frac{(\beta-\alpha)_l}{(\beta)_l}
\,{}_{p+2}F_{q+1} \left[ \left. \begin{gathered}
{-l \quad \alpha \quad (a_p)} \\
{\alpha-\beta-l+1 \quad (b_p)}
\end{gathered} \right| x \right]
\end{equation*}
при   $ p=1, \,\, q=1, \,\, a_1=-l, \,\, b_1=\sigma+\dfrac12, \,\, \alpha=\rho +l+\sigma, \,\, \beta=\alpha+1 $, представим гипергеометрическую функцию Гаусса в \eqref{1_b4} как
\begin{multline*}
{}_2F_1 \left[ \left. \begin{gathered}
{-l \quad \rho+l+\sigma} \\
{\sigma+\dfrac12}
\end{gathered} \right| z^2\right] 
=\frac{\Gamma(\rho+2l+\sigma+1)}{\Gamma(\rho+l+\sigma)}
\times
\\
\times
\sum_{m=0}^{l} \frac{(-1)^m \, \Gamma(\rho+l+m+\sigma)}{m! \, (l-m)! \, \Gamma(\rho+l+m+\sigma+1)} 
{}_{2}F_1 \left[ \left. \begin{gathered}
{-m \quad -l} \\
{\sigma+\dfrac12}
\end{gathered} \right| z^2 \right]
\end{multline*}
и таким образом
\begin{multline}
\label{1_b5}
C_{2l+\sigma}^{\rho}(z) =  
\frac{\Gamma(\rho+2l+\sigma+1)}{\Gamma(\rho)}
\times
\\
\times
\sum_{m=0}^{l} \frac{(-1)^{m+l} \, \Gamma(\rho+l+m+\sigma) (2z)^{\sigma} }{m! \, l! \, (l-m)! \, \Gamma(\rho+l+m+\sigma+1)}
{}_{2}F_1 \left[ \left. \begin{gathered}
{-m \quad -l} \\
{\sigma+\dfrac12}
\end{gathered} \right| z^2 \right]
\end{multline}

С другой стороны покажем, что для $ M $ -- мерных единичных векторов $ \boldsymbol{\zeta}_{1} $ и $ \boldsymbol{\zeta}_{2} $
\begin{multline}
\label{1_b6}
\int\frac{d\Omega_{\boldsymbol{\zeta}}}{S_{_M}} \frac{\left(2(\boldsymbol{\zeta}_{1} \boldsymbol{\zeta})\right)^{2l+\sigma}}{(2l+\sigma)!} \frac{\left(2(\boldsymbol{\zeta}_{2} \boldsymbol{\zeta})\right)^{2m+\sigma}}{(2m+\sigma)!}= \\
=\frac{\Gamma\left(\dfrac{M}{2}\right)}{\Gamma\left(\dfrac{M}{2}+l+m+\sigma\right)} \frac{\left(2(\boldsymbol{\zeta}_{1} \boldsymbol{\zeta}_{2})\right)^{\sigma}}{l! \, m!} \,
{}_{2}F_1 \left[ \left. \begin{gathered}
{-m \quad -l} \\
{\sigma+\dfrac12}
\end{gathered} \right| (\boldsymbol{\zeta}_{1} \boldsymbol{\zeta}_{2})^2 \right]
\end{multline}
Действительно, используя выражение из \cite[Гл.5 п.5.3.2(2)]{Prudnikov_3}
\[
\sum_{k=0}^{n} \frac{(-1)^k \, n! \, (2k+v) \, (v)_k}{k! \, (n-k)! \, (n+v+1)_k}
{}_{p+2}F_q \left[ \left. \begin{gathered}
{-k \quad v+k \quad (a_p)} \\
{(b_q)}
\end{gathered}  \right| x \right]=
(v)_{n+1}\frac{\prod\limits_p (a_p )_n}{\prod\limits_q (b_q)_n}x^n
\]
при  $ p=0, \,\, q=1, \,\, b_1=\sigma+\dfrac12, \,\, v=\dfrac{M}{2}+\sigma-1 $  
и используя представление гипергеометрической функции через полиномы Гегенбауэра из \eqref{1_b4} получим, что
\begin{equation}
\label{1_b7}
\begin{gathered}
\sum_{k=0}^{n} \frac{\left(2k+\sigma+\dfrac{M}{2}-1\right)\Gamma\left(\dfrac{M}{2}-1\right)}{(n-k)! \, \Gamma\left(\dfrac{M}{2}+n+k+\sigma\right)}
C_{2k+\sigma}^{\tfrac{M}{2}-1}(z)=\frac{(2z)^{2n+\sigma}}{(2n+\sigma)!} \\
\sum_{k=0}^{n-2k\geqslant 0} \frac{n-2k+\dfrac{M}{2}-1}{\left(\dfrac{M}{2}-1\right) \left(\dfrac{M}{2}\right)_{n-k} k!} C_{n-2k}^{\tfrac{M}{2}-1}(z)=\frac{(2z)^n}{n!}
\end{gathered}
\end{equation}
Таким образом представляя левую часть выражения в \eqref{1_b6}  при $ z=(\boldsymbol{\zeta}_{1} \boldsymbol{\zeta}) $ и $ z=(\boldsymbol{\zeta}_{2} \boldsymbol{\zeta}) $ запишем ее как
\begin{multline*}
\int\frac{d\Omega_{\boldsymbol{\zeta}}}{S_{_M}} \frac{\left(2(\boldsymbol{\zeta}_{1} \boldsymbol{\zeta})\right)^{2l+\sigma}}{(2l+\sigma)!} \frac{\left(2(\boldsymbol{\zeta}_{2} \boldsymbol{\zeta})\right)^{2m+\sigma}}{(2m+\sigma)!}= \\
=\sum_{k_1=0}^{l} \sum_{k_2=0}^{m} \frac{\left(2k_1+\sigma+\dfrac{M}{2}-1\right) \left(2k_2+\sigma+\dfrac{M}{2}-1\right) \Gamma\left(\dfrac{M}{2}-1\right)^2}{(l-k_1)! \, (m-k_2)! \, \Gamma\left(\dfrac{M}{2}+l+k_1+\sigma\right) \Gamma\left(\dfrac{M}{2}+m+k_2+\sigma\right)} \times \\
\times \int\frac{d\Omega_{\boldsymbol{\zeta}}}{S_{_M}}
C_{2k_1+\sigma}^{\tfrac{M}{2}-1}\left((\boldsymbol{\zeta}_{1} \boldsymbol{\zeta})\right) C_{2k_2+\sigma}^{\tfrac{M}{2}-1}\left((\boldsymbol{\zeta}_{2} \boldsymbol{\zeta})\right)
\end{multline*}
Используя теорему о свертке для сферических гармоник \cite[Гл.11]{Beitman_2}, в котором
\begin{equation*}
\int\frac{d\Omega_{\boldsymbol{\zeta}}}{S_{_M}}
C_{2k_1+\sigma}^{\tfrac{M}{2}-1}\left((\boldsymbol{\zeta}_{1} \boldsymbol{\zeta})\right) C_{2k_2+\sigma}^{\tfrac{M}{2}-1}\left((\boldsymbol{\zeta}_{2} \boldsymbol{\zeta})\right)=
\delta_{k_1,k_2}\frac{\dfrac{M}{2}-1}{2k_1+\sigma+\dfrac{M}{2}-1} C_{2k_1+\sigma}^{\tfrac{M}{2}-1}\left((\boldsymbol{\zeta}_{1} \boldsymbol{\zeta}_{2})\right)
\end{equation*}
представим правую часть предыдущего выражение как
\begin{multline*}
\sum_{k_1=0}^{\infty} \frac{\left(2k_1+\sigma+\dfrac{M}{2}-1\right) \Gamma\left(\dfrac{M}{2}-1\right)^2 \left(\dfrac{M}{2}-1\right) }{ \Gamma(l+1-k_1) \Gamma(m+1-k_1) \Gamma\left(\dfrac{M}{2}+l+k_1+\sigma\right) \Gamma\left(\dfrac{M}{2}+m+k_1+\sigma\right)} 
\times 
\\ 
\times
C_{2k_1+\sigma}^{\tfrac{M}{2}-1}\left((\boldsymbol{\zeta}_{1} \boldsymbol{\zeta}_{2})\right)
\end{multline*}
выражая в этом выражении функцию Гегенбауэра в виде ряда \eqref{1_b4} (при $ \rho \to \dfrac{M}{2}-1 $  )
\[
C_{2k_1+\sigma}^{\tfrac{M}{2}-1}\left((\boldsymbol{\zeta}_{1} \boldsymbol{\zeta}_{2})\right)=
\sum_{n=0}^{\infty} \frac{(-1)^{k_1+n}\Gamma \left( \dfrac{M}{2}-1+k_1+n+\sigma \right) )}{\Gamma(k_1+1-n)\Gamma\left(\dfrac{M}{2}-1 \right) } \frac{\left( 2(\boldsymbol{\zeta}_{1} \boldsymbol{\zeta}_{2})\right)^{2n+\sigma}}{(2n+\sigma)!}
\]
запишем предпоследнее выражение в виде
\begin{multline*}
\sum_{n=0}^{\infty} \sum_{k_1=0}^{\infty}
\frac{(-1)^{k_1+n} \left(2k_1+\sigma+\dfrac{M}{2}-1\right) \Gamma\left(\dfrac{M}{2}\right) }{\Gamma(l+1-k_1) \Gamma(m+1-k_1) \Gamma(k_1+1-n)  } 
\times 
\\
\times 
\frac{\Gamma\left(\dfrac{M}{2}-1+k_1+n+\sigma \right)}{\Gamma\left(\dfrac{M}{2}+l+k_1+\sigma\right) \Gamma\left(\dfrac{M}{2}+m+k_1+\sigma\right)} 
\frac{\left( 2(\boldsymbol{\zeta}_{1} \boldsymbol{\zeta}_{2})\right)^{2n+\sigma}}{(2n+\sigma)!} 
\end{multline*}
после суммирования по $ k_1 $   с учетом того, что $ n $   принимает целые положительные числа, получим
\begin{multline*}
\sum_{n=0}^{\infty}\frac{\Gamma\left(\dfrac{M}{2}\right) \Gamma(a+1)}{\Gamma(l+1-n)\Gamma(m+1-n)\Gamma(1+a-b)\Gamma(1+a-c)} \times \\ \times
{}_4F_3 \left[ \left. \begin{gathered}
{a \quad \dfrac{a}{2}+1 \quad b \quad c} \\
{\dfrac{a}{2} \quad 1+a-b \quad 1+a-c }
\end{gathered} \right| -1 \right]
\frac{\left( 2(\boldsymbol{\zeta}_{1} \boldsymbol{\zeta}_{2})\right)^{2n+\sigma}}{(2n+\sigma)!}
\end{multline*}
\begin{equation*}
\text{где} \quad a=-1+2n+\sigma+\dfrac{M}{2}, \quad b=n-m, \quad c=n-l
\end{equation*}
гипергеометрический ряд $ {}_4F_3 \left[ \left. \cdots \right| -1 \right] $ можно выразить из \cite[Гл.7  п.7.5.4(2)]{Prudnikov_3} как
\[
{}_4F_3 \left[ \left. \begin{gathered}
{a \quad \dfrac{a}{2}+1 \quad b \quad c} \\
{\dfrac{a}{2} \quad 1+a-b \quad 1+a-c }
\end{gathered} \right| -1 \right]=\frac{\Gamma(1+a-b)\Gamma(1+a-c)}{\Gamma(1+a)\Gamma(1+a-b-c)}, \quad a-2b-2c>-1
\]
и окончательно получим
\begin{gather*}
\sum_{n=0}^{\infty} \frac{\Gamma\left(\dfrac{M}{2}\right)}{\Gamma(l+1-n)\Gamma(m+1-n)\Gamma\left(\dfrac{M}{2}+l+m+\sigma\right)}
\frac{\left( 2(\boldsymbol{\zeta}_{1} \boldsymbol{\zeta}_{2})\right)^{2n+\sigma}}{(2n+\sigma)!} =\\
=\frac{\Gamma\left(\dfrac{M}{2}\right)}{\Gamma\left(\dfrac{M}{2}+l+m+\sigma\right)} 
\frac{\left( 2(\boldsymbol{\zeta}_{1} \boldsymbol{\zeta}_{2})\right)^{\sigma}}{l! \, m!}
{}_{2}F_1 \left[ \left. \begin{gathered}
{-m \quad -l} \\
{\sigma+\dfrac12}
\end{gathered} \right| (\boldsymbol{\zeta}_{1} \boldsymbol{\zeta}_{2})^2 \right] \\
\dfrac{M}{2}+2l+2m-2n+\sigma=\dfrac{M}{2}+l+m+\sigma+(l-n)+(m-n)>0  
\end{gather*}
что и доказывает выражение \eqref{1_b6}. Объединяя \eqref{1_b5} при $ z{=}(\boldsymbol{\zeta}_{1} \boldsymbol{\zeta}_{2}) $   и \eqref{1_b6} получим \eqref{1_b3}.

Используя \eqref{1_b3} запишем \eqref{1_b2} как
\begin{multline*}
\frac{1}{\left|\mathbf{R}_{_{N-1}}+\mathbf{r}_{_N}\right|^{2\rho}}=
\frac{1}{(r_{_N})^{2\rho}} \int\frac{d\Omega_{\boldsymbol{\zeta}}}{S_{_M}} 
\sum_{\sigma=0}^{1}\sum_{l=0}^{\infty}\frac{\Gamma(\rho+2l+\sigma+1)}{(2l+\sigma)! \, \Gamma(\rho)}
\left( 2\frac{R_{_{N-1}}}{r_{_N}} (\boldsymbol{\zeta}_{R} \boldsymbol{\zeta}) \right)^{2l+\sigma} \times \\ \times
\sum_{m=0}^{l}
\frac{(-1)^{l+m+\sigma} \, \Gamma(\rho+l+m+\sigma) \, \Gamma\left(\dfrac{M}{2}+l+m+\sigma\right)}{(l-m)! \, \Gamma(\rho+l+m+\sigma+1) \, \Gamma\left(\dfrac{M}{2}\right)}
\frac{\left( 2(\boldsymbol{\zeta}_{N} \boldsymbol{\zeta})\right)^{2m+\sigma}}{(2m+\sigma)!}
\end{multline*}

Объединяя четные ( $ \sigma=0 $ ) и нечетные ($ \sigma=1 $ ) суммы по $ n=2l+\sigma=0,1,2,\ldots $, преобразуем в этом выражений суммы как
\begin{multline}
\label{1_b9}
\frac{1}{\left|\mathbf{R}_{_{N-1}}+\mathbf{r}_{_N}\right|^{2\rho}}=
\frac{1}{(r_{_N})^{2\rho}} \int\frac{d\Omega_{\boldsymbol{\zeta}}}{S_{_M}} 
\sum_{n=0}^{\infty}\sum_{m=0}^{n-2m\geqslant 0}\frac{\Gamma(\rho+n+1)}{n! \, \Gamma(\rho)}
\left( 2\frac{R_{_{N-1}}}{r_{_N}} (\boldsymbol{\zeta}_{R} \boldsymbol{\zeta}) \right)^{n} 
\times 
\\ 
\times
\frac{(-1)^{n+m} \, \Gamma(\rho+n-m) \, \Gamma\left(\dfrac{M}{2}+n-m\right)}{m! \, \Gamma(\rho+n-m+1) \, \Gamma\left(\dfrac{M}{2}\right)}
\frac{\left( 2(\boldsymbol{\zeta}_{N} \boldsymbol{\zeta})\right)^{n-2m}}{(n-2m)!}
\end{multline} 
Используя соотношение вида
\begin{multline*}
\sum_{n=0}^{\infty}\frac{\prod\limits_p \Gamma(a_p+n)}{\prod\limits_q \Gamma(b_q+n)} \frac{\alpha^n}{n!}
\left(\sum_{i=1}^{N-1} z_i\right)^n = \\
= \sum_{n_1,\ldots,n_{_{N-1}}=0}^{\infty} \,\,\, 
\alpha^{n_1+n_2+\ldots +n_{_{N-1}}}
\frac{\prod\limits_p \Gamma(a_p+n_1+n_2+\ldots +n_{_{N-1}})}{\prod\limits_q \Gamma(b_q+n_1+n_2+\ldots +n_{_{N-1}})}
\prod_{i=1}^{N-1}\frac{(z_i)^{n_i}}{n_i !}
\end{multline*}
где для сокращения записи мы введем 
\begin{gather*}
\tilde{n}=n_1+\ldots +n_{_{N-1}} =\sum\limits_{i=1}^{N-1}n_i , \quad 
\sum\limits_{n_1=0}^{\infty} \,\sum\limits_{n_2=0}^{\infty}\ldots \sum\limits_{n_{_{N-1}}=0}^{\infty} (\ldots) =\sum\limits_{\tilde{n}}(\ldots) 
\end{gather*}
Выражая в \eqref{1_b9} при
$ \mathbf{R}_{_{N-1}}{=}R_{_{N-1}} \boldsymbol{\zeta}_{R} {=} \mathbf{r}_1{+}\mathbf{r}_2{+}\ldots{+}\mathbf{r}_{_{N-1}}{=}\sum\limits_{i=1}^{N-1}r_i \boldsymbol{\zeta}_{i} $   представим \\  
$ 2\dfrac{R_{_{N-1}}}{r_{_N}} (\boldsymbol{\zeta}_{R} \boldsymbol{\zeta}) =\sum\limits_{i=1}^{N-1}2\dfrac{r_i}{r_{_N}} (\boldsymbol{\zeta}_{i} \boldsymbol{\zeta})=\sum\limits_{i=1}^{N-1} z_i $.
Тогда используя последнее соотношение запишем  \eqref{1_b9}  как 
\begin{multline}
\label{1_b10}
\frac{1}{\left| \left( \mathbf{r}_1+\ldots +\mathbf{r}_{_{N-1}} \right) +\mathbf{r}_{_N}\right|^{2\rho}}=
\\
=\frac{1}{(r_{_N})^{2\rho}} \int\frac{d\Omega_{\boldsymbol{\zeta}}}{S_{_M}} 
\sum_{\tilde{n}}\prod_{i=1}^{N-1} \left[
\left( \frac{r_i}{r_{_N}}\right)^{n_i} \frac{\left(2(\boldsymbol{\zeta}_{i} \boldsymbol{\zeta})\right)^{n_i}}{n_i !}
\right]
\frac{(-1)^{\tilde{n}}\Gamma(\rho+\tilde{n}+1)}{\Gamma(\rho)} 
\times 
\\ 
\times 
\sum_{m=0}^{\tilde{n}-2m\geqslant 0}\frac{(-1)^{m} \, \Gamma(\rho+\tilde{n}-m) \, \Gamma\left(\dfrac{M}{2}+\tilde{n}-m\right)}{m! \, \Gamma(\rho+\tilde{n}-m+1) \, \Gamma\left(\dfrac{M}{2}\right)}
\frac{\left( 2(\boldsymbol{\zeta}_{N} \boldsymbol{\zeta})\right)^{\tilde{n}-2m}}{(\tilde{n}-2m)!}
\end{multline}
Рассмотрим отдельную сумму по $ m $ . Так, используя второе выражение в \eqref{1_b7} \\ при  $ z=(\boldsymbol{\zeta}_{N} \boldsymbol{\zeta}) $ представим его как
\begin{multline*}
\sum_{m=0}^{\tilde{n}{-}2m\geqslant 0} \, \sum_{k=0}^{\tilde{n}{-}2m{-}2k\geqslant 0}
\frac{({-}1)^m \, \Gamma(\rho{+}\tilde{n}{-}m) \Gamma\left(\dfrac{M}{2}{+}\tilde{n}{-}m\right) \left(\dfrac{M}{2}{+}\tilde{n}{-}2m{-}2k{-}1\right) }{m! \, k! \, \Gamma(\rho{+}\tilde{n}{-}m{+}1) \Gamma\left(\dfrac{M}{2}{+}\tilde{n}{-}2m{-}k\right) \left(\dfrac{M}{2}{-}1\right) } \times \\ \times
C_{\tilde{n}{-}2m{-}2k}^{\tfrac{M}{2}{-}1} \left((\boldsymbol{\zeta}_{N} \boldsymbol{\zeta})\right)=
\end{multline*}
\begin{multline*}
=\sum_{\mu_{_N}=0}^{\tilde{n}{-}2\mu_{_N}\geqslant 0} \, \sum_{k=0}^{\mu_{_N}}
\frac{({-}1)^{k{+}\mu_{_N}} \, \Gamma(\rho{+}\tilde{n}{-}\mu_{_N}{+}k) \Gamma\left(\dfrac{M}{2}{+}\tilde{n}{-}\mu_{_N}{+}k\right) \left(\dfrac{M}{2}{+}\tilde{n}{-}2\mu_{_N}{-}1\right) }{(\mu_{_N}{-}k)! \, k! \, \Gamma(\rho{+}\tilde{n}{-}\mu_{_N}{+}k{+}1) \Gamma\left(\dfrac{M}{2}{+}\tilde{n}{-}2\mu_{_N}{+}k\right) \left(\dfrac{M}{2}{-}1\right) } \times \\ \times
C_{\tilde{n}{-}2\mu_{_N}}^{\tfrac{M}{2}{-}1} \left((\boldsymbol{\zeta}_{N} \boldsymbol{\zeta})\right)
\end{multline*}
где во втором выражении мы преобразовали порядок суммирования. Снимая сумму по $ k $  получим
\begin{multline*}
\sum_{\mu_{_N}=0}^{\tilde{n}-2\mu_{_N}\geqslant 0} \frac{(-1)^{\mu_{_N}}}{\mu_{_N} !}
{}_3F_2 \left[ \left. \begin{gathered}
{-\mu_{_N} \quad \rho+\tilde{n}-\mu_{_N} \quad \dfrac{M}{2}+\tilde{n}-\mu_{_N}} \\
{\rho+\tilde{n}-\mu_{_N}+1 \quad \dfrac{M}{2}+\tilde{n}-2\mu_{_N}}
\end{gathered} \right| 1 \right]
\times \\ \times
\frac{ \Gamma(\rho+\tilde{n}-\mu_{_N}) \Gamma\left(\dfrac{M}{2}+\tilde{n}-\mu_{_N}\right) \left(\dfrac{M}{2}+\tilde{n}-2\mu_{_N}-1\right) }{ \Gamma(\rho+\tilde{n}-\mu_{_N}+1) \Gamma\left(\dfrac{M}{2}+\tilde{n}-2\mu_{_N}\right) \left(\dfrac{M}{2}-1\right) } 
C_{\tilde{n}-2\mu_{_N}}^{\tfrac{M}{2}-1} \left((\boldsymbol{\zeta}_{N} \boldsymbol{\zeta})\right)
\end{multline*}
здесь гипергеометрический ряд $ {}_3F_2 \left[ \left. \ldots \right| 1 \right] $   является рядом  Заальшютса \cite[Гл.4 п.4.4(3)]{Beitman_1}
{\small \begin{equation*}
	{}_3F_2 \left[ \left. \begin{gathered}
	{-\mu_{_N} \quad \rho{+}\tilde{n}{-}\mu_{_N} \quad \dfrac{M}{2}{+}\tilde{n}{-}\mu_{_N}} \\
	{\rho{+}\tilde{n}{-}\mu_{_N}{+}1 \quad \dfrac{M}{2}{+}\tilde{n}{-}2\mu_{_N}}
	\end{gathered} \right| 1 \right]
	=
	\frac{\mu_{_N} ! \, \left(\rho{+}1{-}\dfrac{M}{2}\right)_{\mu_{_N}}}{(\rho{+}\tilde{n}{-}\mu_{_N}{+}1)_{\mu_{_N}} \left({-}\dfrac{M}{2}{-}\tilde{n}{+}\mu_{_N}{+}1 \right)_{\mu_{_N}}}
	\end{equation*}}
и таким образом сумму по $ m $  в \eqref{1_b10} можно представить в виде
\begin{equation*}
\sum_{\mu_{_N}=0}^{\tilde{n}-2\mu_{_N}\geqslant 0}
\frac{ \Gamma(\rho+\tilde{n}-\mu_{_N}) \Gamma\left(\rho-\dfrac{M}{2}+\mu_{_N}+1\right) }{ \Gamma(\rho+\tilde{n}+1) \Gamma\left(\rho-\dfrac{M}{2}+1\right)  }
\frac{\dfrac{M}{2}+\tilde{n}-2\mu_{_N}-1}
{\dfrac{M}{2}-1} 
C_{\tilde{n}-2\mu_{_N}}^{\tfrac{M}{2}-1} \left((\boldsymbol{\zeta}_{N} \boldsymbol{\zeta})\right)
\end{equation*}
Выражая в \eqref{1_b10} $ \dfrac{\left(2(\boldsymbol{\zeta}_{i} \boldsymbol{\zeta})\right)^{n_i}}{n_i !} $   из второго выражения \eqref{1_b7} при $z=(\boldsymbol{\zeta}_{i} \boldsymbol{\zeta})  $   и объединяя с последним выражением (сумма по $ m $ в \eqref{1_b10}) получим
\begin{multline*}
\frac{1}{\left|\mathbf{r}_1+\ldots +\mathbf{r}_{_{N-1}}+\mathbf{r}_{_N}\right|^{2\rho}}=
\\
=\frac{1}{(r_{_N})^{2\rho}} \int\frac{d\Omega_{\boldsymbol{\zeta}}}{S_{_M}} 
\sum_{\tilde{n}} (-1)^{\tilde{n}} (\rho)_{\tilde{n}-\mu_{_N}} \left( \rho-\dfrac{M}{2}+1\right)_{\mu_{_N}}
\sum_{\mu_1=0}^{\left[\frac{n_1}{2}\right]} \ldots \sum_{\mu_{_{N-1}}=0}^{\left[\frac{n_{_{N-1}}}{2}\right]} \,\, \sum_{\mu_{_N}=0}^{\left[\frac{\tilde{n}}{2}\right]}
\times 
\\ 
\times
\prod_{i=1}^{N-1} 
\left(
\frac{ \left(\dfrac{r_i}{r_{_N}}\right)^{n_i}}{\mu_i ! \, \left(\dfrac{M}{2}\right)_{n_i-\mu_i}} \frac{n_i-2\mu_i+\dfrac{M}{2}-1}{\dfrac{M}{2}-1}C_{n_i-2\mu_i}^{\tfrac{M}{2}-1}\left((\boldsymbol{\zeta}_{i} \boldsymbol{\zeta})\right) 
\right)
\times \\ \times
\frac{\tilde{n}-2\mu_{_N}+\dfrac{M}{2}-1}{\dfrac{M}{2}-1}
C_{\tilde{n}-2\mu_{_N}}^{\tfrac{M}{2}-1}\left((\boldsymbol{\zeta}_{N} \boldsymbol{\zeta})\right) 
\end{multline*}
или вводя $ n_i-2\mu_i=l_i, \, i=1,\ldots , N{-}1 $   и  $ \tilde{n}-2\mu_{_N}=l_{_N} $ а также $ l{=}\dfrac12 \sum\limits_{i=1}^{N}l_i, \,\, \nu {=}2\rho $, преобразуем суммы в последнем выражении и запишем в более простом виде как в \eqref{1_b12}
\begin{multline*}
\frac{1}{\left|\mathbf{r}_1+\ldots+\mathbf{r}_{_N}\right|^{\nu}}= \\
=\frac{1}{(r_{_N})^{\nu}} \sum_{\substack{l_1=0 \\\dots \\ l_{_N}=0}}^{\infty}  \,\, \sum_{\substack{\mu_1=0 \\ \dots \\ \mu_{_{N-1}}=0}}^{\infty}
\frac{(-1)^{l_{_N}} \left(\dfrac{\nu}{2}\right)_{\mu_1+\ldots +\mu_{_{N-1}}+l} \left(\dfrac{\nu{-}M{+}2}{2}\right)_{\mu_1+\ldots +\mu_{_{N-1}}+l-l_{_N}} }{\displaystyle\prod\limits_{i=1}^{N-1} \left[\left(\dfrac{M}{2}\right)_{l_i+\mu_i} \mu_i ! \right]}
\times \\ \times
\prod\limits_{i=1}^{N-1} \left(\frac{r_i}{r_{_N}}\right)^{l_i+2\mu_i}
\int\frac{d\Omega_{\boldsymbol{\zeta}}}{S_{_M}} \prod\limits_{i=1}^{N} \frac{l_i+\dfrac{M}{2}-1}{\dfrac{M}{2}-1} C_{l_i}^{\tfrac{M}{2}-1}\left((\boldsymbol{\zeta}_{i} \boldsymbol{\zeta})\right)
\end{multline*}
Используя теорему сложения \cite[Гл.11, п.11.4]{Beitman_2} 
\begin{equation} \label{th_C}
C_l^{\tfrac{M}{2}-1} \left( ( \boldsymbol{\zeta}_{1} \boldsymbol{\zeta}_{2} ) \right) = 
S_{_M} \frac{\dfrac{M}{2}-1}{l+\dfrac{M}{2}-1} \sum_{\mathbf{m}_{\mathbf{k}}} 
Y_{l,\mathbf{m}_{\mathbf{k}}} \left(\boldsymbol{\zeta}_{1} \right) 
Y_{l,\mathbf{m}_{\mathbf{k}}} \left(\boldsymbol{\zeta}_{2} \right)^{\ast}
\end{equation}
получим коэффициентную угловую функцию вида \eqref{1_b11}.

Выражение \eqref{1_b12} можно так же представить через функцию Лауричеля многих переменных как
\begin{equation}\label{1_bR}
\frac{1}{\left|\mathbf{r}_1+\ldots+\mathbf{r}_{_N}\right|^{\nu}}= \sum_{l_1,\ldots,l_{_N}=0}^{\infty} 
V_{l_1,\ldots, l_{_N}} \left( \boldsymbol{\zeta}_{1},\ldots ,\boldsymbol{\zeta}_{N} \right)
R_{l_1,\ldots, l_{_N}}^{(\nu,M)}(r_1,\ldots,r_{_{N}})
\end{equation}
где
\begin{multline} \label{1_bL}
R_{l_1,\ldots, l_{_N}}^{(\nu,M)}(r_1,\ldots,r_{_{N}})=
\frac{(-1)^{l_{_N}}}{(r_{_N})^{\nu}}  
\frac{\left(\dfrac{\nu}{2}\right)_{l} \left(\dfrac{\nu-M+2}{2}\right)_{l-l_{_N}} }{\displaystyle\prod\limits_{p=1}^{N-1} \left(\dfrac{M}{2}\right)_{l_p}}
\times \\ \times
\prod\limits_{p=1}^{N-1} \left(\frac{r_p}{r_{_N}}\right)^{l_p}
F_C^{(N-1)} \left[ \left. \begin{gathered}
{l+\dfrac{\nu}{2} \quad \dfrac{\nu -M+2}{2}+l-l_{_N}} \\
{l_1+\dfrac{M}{2},\ldots ,l_{_{N-1}}+\dfrac{M}{2}}
\end{gathered} \right|\left(\frac{r_1}{r_{_N}}\right)^{2},\ldots ,\left(\frac{r_{_{N-1}}}{r_{_N}}\right)^{2}  \right] 
\end{multline}
Сходимость этого ряда и ряда \eqref{1_b12}, такая же как и в \eqref{1_b2} выполняется при условии $ \sum\limits_{p=1}^{N-1}\dfrac{r_p}{r_{_N}}<1 $ . Действительно, из \eqref{1_b2} сходимость выполняется при условии, что $ \dfrac{\left|\mathbf{R}_{_{N-1}}\right|}{r_{_N}}<1 $   при $ \max \left|\mathbf{R}_{_{N-1}}\right|= \sum\limits_{p=1}^{N-1} r_p $

\section{Представление вида \eqref{1_a015} для $ N>2 $ многомерных векторов} \label{sec_1_3}

\begin{Theorm} \label{Th_1_3}
	Для произвольных $ M $ ---мерных векторов $ \mathbf{r}_{s}=r_s \cdot \boldsymbol{\zeta}_{s} \quad s=1,2,\ldots,N $, где $ \boldsymbol{\zeta}_{s} $ единичные $ M $ ---мерные вектора в евклидовом пространстве, а также  действительных чисел  $ 0<\nu<M $, для однородных функции  $ \left| \mathbf{r}_1{+}\ldots {+}\mathbf{r}_{_N} \right|^{-\nu} \in \mathbb{R} $ справедливо выражение
	\begin{multline} \label{eq:c20}
	\frac{1}{\left|\mathbf{r}_1+\ldots+\mathbf{r}_{_N}\right|^{\nu}}=\frac{2\Gamma\left(\dfrac{M-\nu}{2}\right) }
	{S_{_M} \, \Gamma\left(\dfrac{M}{2}\right)\Gamma\left(\dfrac{\nu}{2}\right)}  \sum_{\mathbf{n,l},\mathbf{m}_{\mathbf{k}}^{(i)}} \,\, \prod_{p=1}^{N}
	\left( \left(r_p^2+1\right)^{\tfrac{M-\nu}{2}} H_{n_p,l_p,\mathbf{m}_{\mathbf{k}}^{(p)}}^{(M)}(\mathbf{r}_p)\right)
	\times
	\\
	\times
	\int\limits_{0}^{\infty}du \, u^{\nu-1} \int d\Omega_{\boldsymbol{\zeta}} \prod_{s=1}^{N}
	\Xi_{n_s,l_s,\mathbf{m}_{\mathbf{k}}^{(s)}}^{(\nu,M)}(u,\boldsymbol{\zeta})
	\end{multline}
	где 
	\begin{equation} \label{eq:c19}
	\Xi_{n,l,\mathbf{m}_{\mathbf{k}}}^{(\nu,M)}(u,\boldsymbol{\zeta})=
	\xi_{n,l}^{(\nu,M)}(u) Y_{l,\mathbf{m}_{\mathbf{k}}}(\boldsymbol{\zeta})^{\ast}
	\end{equation}
	функция 
	\begin{multline} \label{eq:c16}
	\xi_{n,l}^{(\nu,M)}(u)= 
	\frac{i^{l}2\pi^{\tfrac{M}{2}} u^{l}}
	{\Gamma\left(l{+}\dfrac{M}{2}\right) \Gamma\left(l{+}M{-}\dfrac{\nu}{2}\right)}
	\sqrt{\frac{\left(n{+}l{+}\dfrac{M{-}1}{2}\right)\Gamma(n{+}2l{+}M{-}1)}{\Gamma(n{+}1)}}
	\times  
	\\
	\times
	\int\limits_{0}^{\infty}dz \,e^{-z-\tfrac{u^2}{z}} z^{\tfrac{M-\nu}{2}-1} 
	{}_2F_2 \left[ \left. 
	\begin{gathered}
	{ -n \quad n+2l+M-1} \\
	{ l +\dfrac{M}{2} \quad l + M - \dfrac{\nu}{2} } 
	\end{gathered} \right| z  \right]
	\end{multline}
	образует ортогональную систему с весом $ u^{\nu-1} $ в области $ u\in [0\ldots \infty) $
	\begin{equation} \label{eq:c28}
	\int\limits_{0}^{\infty}\!\!du \, u^{\nu-1} 
	\xi_{n_1,l}^{(\nu,M)}(u) \xi_{n_2,l}^{(\nu,M)}(u) =
	\frac{(-1)^l \, \pi^{M} \, \Gamma\left(\dfrac{\nu}{2}+l+n_2\right) }
	{\Gamma\left(M-\dfrac{\nu}{2}+l+n_2\right) }
	\delta_{n_1,n_2}       
	\end{equation}
	И
	\begin{equation} \label{eq:eH}
	H_{n,l,\mathbf{m}_{\mathbf{k}}}^{(M)} (\mathbf{r}_s)=
	\eta_{n,l}^{(M)}(r_s) Y_{l,\mathbf{m}_{\mathbf{k}}} \left(\boldsymbol{\zeta}_{s}\right)
	\end{equation}
	где функция 
	\begin{multline} \label{eq:e_eta}
	\eta_{n,l}^{\left(M\right)}(r)=
	2^{2l+M-1}\Gamma\left( l{+}\frac{M{-}1}{2} \right)\sqrt{\frac{n!\,\left(n{+}l{+}\dfrac{M{-}1}{2}\right)}{\pi\,\Gamma(n{+}2l{+}M{-}1)}} \frac{r^l}{\left(r^2{+}1\right)^{l+\tfrac{M}{2}}}
	C_n^{l+\tfrac{M-1}{2}}\left( \frac{r^2{-}1}{r^2{+}1} \right)=
	\\
	\shoveleft{=\frac{2}{\Gamma\left(l{+}\dfrac{M}{2}\right)}
	\sqrt{\frac{\left(n{+}l{+}\dfrac{M{-}1}{2}\right)\Gamma(n{+}2l{+}M{-}1)}{n!}} 
	\times} 
	\\
	\times \frac{r^l}{\left(r^2{+}1\right)^{l+\tfrac{M}{2}}}
	\,{}_2F_1 \left[ \left. 
	\begin{matrix}
	{  -n  \quad n{+}2l{+}M{-}1  } \\
	{  l{+}\dfrac{M}{2}  }
	\end{matrix} \right| \frac{1}{r^2{+}1}\right] 
	\end{multline}
	образует ортогональную систему с весом $ r^{M-1} $
	\begin{equation} \label{eq:e_int_eta}
	\int\limits_{0}^{\infty}\!\!dr \, r^{M-1} \eta_{n_1,l}^{(M)}(r) \eta_{n_2,l}^{(M)}(r)=\delta_{n_1,n_2}
	\end{equation}
\end{Theorm}

В \eqref{eq:c20} общая сумма означает как (где $ \mathbf{n}=\{n_1,\ldots, n_{_N}\},\,\,\mathbf{l}=\{l_1,\ldots, l_{_N}\} $)
\[ 
\sum_{\mathbf{n,l},\mathbf{m}_{\mathbf{k}}^{(i)}} = 
\sum_{n_1=0}^{\infty} \ldots \sum_{n_{_N}=0}^{\infty} \,\,\,
\sum_{l_1=0}^{\infty} \ldots \sum_{l_{_N}=0}^{\infty} \,\,\,
\sum_{\mathbf{m}_{\mathbf{k}}^{(1)},\ldots ,\mathbf{m}_{\mathbf{k}}^{(N)}}
\]

Отметим, что \eqref{1_eH3} и \eqref{1_e_eta3} соответствуют \eqref{eq:eH} и \eqref{eq:e_eta} при $ M=3 $, а  \eqref{1_a015} является частным случаем теоремы \ref{Th_1_3}  при $ M=3,\,N=2 $. 

Так как в \eqref{eq:c16} есть комплексный множитель $ i^l $, выражение \eqref{eq:c20} будет принимать действительные значения. Действительно, как было показано в разделе \ref{sec_1_2} из \eqref{1_b11}, или после интегрирования в \eqref{eq:c20} по $ d\Omega_{\boldsymbol{\zeta}} $ все суммы по $ l_s $ будут четными, и следовательно $ \prod\limits_{s=1}^{N}i^{l_s}=i^{2l}=(-1)^l $. 

Для вывода \eqref{eq:c20} будем исходить из \eqref{1_bR} и \eqref{1_bL}. Так, используя выражение интеграла от произведений для функции Беселя $ J_{\lambda_k}(c_ku) $ из \cite[Гл.2, п.2.12.44(7)]{Prudnikov_3}  ($ J_{\lambda}(z) $---функция Бесселя)
\begin{multline} \label{eq:eJ} 
\int\limits_{0}^{\infty}du \, u^{\alpha-1}
\prod_{k=1}^{N} J_{\lambda_k}(c_k u)=
\frac{2^{\alpha-1}\left( c_{_N} \right)^{\lambda_{_N}-\beta} \Gamma\left( \dfrac{\beta}{2} \right)}{\Gamma\left( \lambda_{_N}-\dfrac{\beta}{2}+1 \right)}
\times
\\
\times
\prod_{k=1}^{N-1}\frac{(c_k)^{\lambda_k}}{\Gamma( \lambda_k +1 )}
\,F_{C}^{(N-1)} \left[ 
\left.
\begin{gathered}
{ \frac{\beta}{2} \quad \frac{\beta}{2}-\lambda_{_N}  } \\
{ \lambda_1+1 ,\ldots , \lambda_{_{N-1}}+1 }
\end{gathered}
\right|  \left( \frac{c_1}{c_{_N}} \right)^2,\ldots,\left( \frac{c_{_{N-1}}}{c_{_N}} \right)^2  \right]
\end{multline}
где при действительных $ \alpha,\,\lambda_k,\,c_k $
\begin{gather*}
N\geqslant 2; \quad c_k >0, \,\, k=1,\ldots,N; \quad c_{_N}>c_1+ \ldots + c_{_{N-1}}; \\
-(\lambda_1+\ldots+\lambda_{_N})<\alpha <\frac{N}{2}+1; \quad \beta=\alpha+\lambda_1+\ldots+\lambda_{_N}
\end{gather*}
Сопоставляя это выражение с \eqref{1_bL} при значениях
\begin{gather*}
c_k=\gamma r_k, \,\, \lambda_k=l_k+\frac{M}{2}-1, \,\, k=1,\ldots,N; \\
\alpha=\nu-N\left( \frac{M}{2}-1 \right)=\nu-\sum_{k=1}^{N}\left( \frac{M}{2}-1 \right); \\ 
\left( \beta=l_1+\ldots+l_{_N}+\nu=2l+\nu; \,\, 
0<\nu<\frac{N(M-1)}{2}+1 \right)
\end{gather*}
получим
\begin{multline} \label{eq:e4}
R_{l_1,\ldots,l_{_N}}^{(\nu,M)} (r_1,\ldots ,r_{_N})=
\\ 
=\frac{(-1)^{l}\,2\,\Gamma\left( \dfrac{M-\nu}{2} \right)}{\Gamma\left( \dfrac{\nu}{2} \right)\Gamma\left( \dfrac{M}{2} \right)} \left( \frac{\gamma}{2} \right)^{\nu}
\int\limits_{0}^{\infty}du\,u^{\nu-1}\prod_{k=1}^{N}\frac{\Gamma\left( \dfrac{M}{2} \right)J_{l_k+\frac{M}{2}-1}(\gamma r_k u)}{\left( \dfrac{\gamma r_k u}{2} \right)^{\frac{M}{2}-1}}=
\\
=\frac{2\,\Gamma\left( \dfrac{M-\nu}{2} \right)}{(S_{_M})^N \Gamma\left( \dfrac{\nu}{2} \right)\Gamma\left( \dfrac{M}{2} \right)} \left( \frac{\gamma}{2} \right)^{\nu}
\int\limits_{0}^{\infty}du\,u^{\nu-1}\prod_{k=1}^{N}\frac{i^{l_k}\, 2\,\pi^{\frac{M}{2}} J_{l_k+\frac{M}{2}-1}(\gamma r_k u)}{\left( \dfrac{\gamma r_k u}{2} \right)^{\frac{M}{2}-1}}  
\end{multline}
где в последнем выражении мы ввели множитель $ S_{_M} $ -- площадь единичной сферы \\ в $ M $--мерном пространстве.
В отличии от \eqref{1_bL}, выражение \eqref{eq:e4} симметрично \\ по всем $ r_k, \, k=1,\ldots,N $.
Вводя функцию вида
\begin{equation}
\label{eq:e5}
\varphi_l^{(M)}(x)=
\frac{i^{l}\, 2\,\pi^{\frac{M}{2}} J_{l+\frac{M}{2}-1}(2x)}{\left(x \right)^{\frac{M}{2}-1}}  
\end{equation}
а также в виду однородности функции \eqref{1_bR}, где $ \gamma $ может быть произвольным (мы положим равным двум) представим \eqref{eq:e4} в следующем виде
\begin{equation}
\label{eq:e6}
R_{l_1,\ldots,l_{_N}}^{(\nu,M)} (r_1,\ldots ,r_{_N})
=\frac{2\,\Gamma\left( \dfrac{M-\nu}{2} \right)}{(S_{_M})^N \Gamma\left( \dfrac{\nu}{2} \right)\Gamma\left( \dfrac{M}{2} \right)}
\int\limits_{0}^{\infty}du\,u^{\nu-1}\prod_{k=1}^{N} 
\varphi_{l_k}^{(M)}(r_k u)  
\end{equation}

\subsection{Ортогонализация $ R_{l_1,\ldots,l_{_N}}^{(\nu,M)}(r_1,\ldots,r_{_N}) $} \label{sec_1_3_1}
Пусть $ \phi_n(r) $ произвольная ортогональная система функций в области значении $ r \in E $ с весом $ \rho(r)>0 $.
\begin{equation} \label{eq:e7}
\int_E dr \, \rho(r)\phi_n(r)\phi_m(r)=\delta_{n,m}
\end{equation}
Представим \eqref{eq:e6} 
\begin{equation}
\label{eq:e8}
R_{\mathbf{l}} (r_1,\ldots ,r_{_N})=\sum_{\mathbf{n}} A_{\mathbf{l,n}}\prod_{k=1}^{N}\phi_{n_k}^{(M)}(r_k)
\end{equation}
где $ \mathbf{l}=\{l_1,l_2,\ldots,l_{_N} \},\,\, \mathbf{n}=\{n_1,n_2,\ldots,n_{_N} \} $
Из \eqref{eq:e7} и \eqref{eq:e8} следует
\begin{equation*}
A_{\mathbf{n,l}}=
\int_E dr_1 \rho(r_1)\phi_{n_1}^{(M)}(r_1) 
\ldots 
\int_E dr_{_N} \rho(r_{_N})\phi_{n_{_N}}^{(M)}(r_{_N})
R_{\mathbf{l}} (r_1,\ldots ,r_{_N})
\end{equation*}
Из \eqref{eq:e6} если ввести функцию вида
\begin{equation*}
\xi_{n,l}^{(\nu,M)}(u)=\int_E dr \rho(r) \, \phi_{n}(r) \varphi_{l}^{(M)}(r u)
\end{equation*}
и используя \eqref{eq:e5} получим
\begin{equation} \label{eq:e9}
\xi_{n,l}^{(\nu,M)}(u)=
i^{l}\, 2\,\pi^{\frac{M}{2}}\int_E dr \rho(r) \phi_{n}(r) \, \frac{ J_{l+\frac{M}{2}-1}(2ru)}{\left(r u \right)^{\frac{M}{2}-1}} 
\end{equation}
\begin{equation}\label{eq:e10}
A_{\mathbf{n,l}}
=\frac{2\,\Gamma\left( \dfrac{M-\nu}{2} \right)}{(S_{_M})^N \Gamma\left( \dfrac{\nu}{2} \right)\Gamma\left( \dfrac{M}{2} \right)}
\int\limits_{0}^{\infty}du\,u^{\nu-1}
\prod_{k=1}^{N} 
\xi_{n_k,l_k}^{(\nu,M)}(u) 
\end{equation}
Найдем такие функции $ \phi_n(r) $ с областью $ E $, при котором функции $ \xi_{n,l}^{(\nu,M)}(u) $ в \eqref{eq:e9} также будут образовывать ортогональную систему функции с весом $ u^{\nu-1} $. Действительно 
\begin{multline} \label{eq:e12}
\int\limits_{0}^{\infty}\!\!du \,u^{\nu-1} 
\xi_{n_1,l}^{(\nu,M)}(u)\xi_{n_2,l}^{(\nu,M)}(u)=
\\
=\int_E\!\! dr_1\,\rho(r_1)\phi_{n_1}(r_1) \int_E\!\! dr_2\,\rho(r_2)\phi_{n_2}(r_2)
\int\limits_{0}^{\infty}\!\!du \,u^{\nu-1}
\varphi_l^{(M)}(r_1 u)\varphi_l^{(M)}(r_2 u)=
\\
=(-1)^l\,4\,\pi^{M}\int_E\!\! dr_1\,\rho(r_1)\phi_{n_1}(r_1) \int_E\!\! dr_2\,\rho(r_2)\phi_{n_2}(r_2)
\times
\\
\times
\int\limits_{0}^{\infty}\!\!du \,u^{\nu-1}
\frac{J_{l+\frac{M}{2}-1}(2r_1u)}{(r_1u)^{\frac{M}{2}-1}}\frac{J_{l+\frac{M}{2}-1}(2r_2u)}{(r_2u)^{\frac{M}{2}-1}}
\end{multline}
Интеграл по $ u $ можно выразить при $ N=2 $ из \eqref{eq:eJ} или найти из общего выражения из \eqref{1_bL} и представить как
\begin{multline} \label{eq:e13}
\int\limits_{0}^{\infty}\!\!du \,u^{\nu-1}
\frac{J_{l+\frac{M}{2}-1}(2r_1u)}{(r_1u)^{\frac{M}{2}-1}}\frac{J_{l+\frac{M}{2}-1}(2r_2u)}{(r_2u)^{\frac{M}{2}-1}}=
\frac{\Gamma\left( l+\dfrac{\nu}{2} \right)}{2\Gamma\left( \dfrac{M-\nu}{2} \right)\Gamma\left( l+\dfrac{M}{2} \right)\,(r_2)^{\nu}}
\times
\\
\times
\left( \frac{r_1}{r_2} \right)^l
\,{}_{2}F_{1} \left[ 
\left.
\begin{gathered}
{ l+\frac{\nu}{2} \quad \frac{\nu-M+2}{2} } \\
{ l+\frac{M}{2} }
\end{gathered}
\right| \left( \frac{r_1}{r_2} \right)^2 \right]=
\\
=\frac{2^{\tfrac{M-\nu-3}{2}}}{\sqrt{\pi}\,\Gamma\left( \dfrac{M-\nu}{2} \right) \left( r_1 r_2 \right)^{\frac{\nu}{2}}}
\frac{e^{-i\,\pi\,\frac{\nu-M+1}{2}}Q_{l+\frac{M-3}{2}}^{\frac{\nu-M+1}{2}}(z)}{\sqrt{z^2-1}^{\frac{\nu-M+1}{2}}}
\end{multline}
\[ z=\frac{r_1^2+r_2^2}{2r_1r_2}>1 \]
где в последнем выражении мы выразили гипергеометрический ряд через функцию Лежандра второго рода \cite[Гл.3, п.3.2(45)]{Beitman_1}.

В работе \cite{AkhR_UNC} было получено (см. Приложение \ref{append}), что для функции Лежандра второго рода справедливо разложение вида
\begin{multline} \label{1_f6}
\frac{e^{-i\,\pi\,\mu}Q_{v}^{\mu}(z)}{\sqrt{z^2-1}^{\mu}}=
\frac{2^{2v+\frac{3}{2}}\Gamma(v+1)^2}{\sqrt{\pi}\left(\sqrt{z_1^2+1}\sqrt{z_2^2+1} \right)^{v+\mu+1}}
\times
\\
\times
\sum_{n=0}^{\infty}
\frac{(-1)^n \, n! \, (n+v+1)}{\Gamma(n+2v+2)}
C_n^{v+1}\left( \frac{z_1}{\sqrt{z_1^2+1}} \right)
C_n^{v+1}\left( \frac{z_2}{\sqrt{z_2^2+1}} \right)
\frac{e^{-i\,\pi\,\left( \mu -\frac{1}{2} \right)}Q_{n+v+\frac{1}{2}}^{\mu-\frac{1}{2}}(z_3)}{\sqrt{z_3^2-1}^{\mu -\frac{1}{2}}}
\end{multline}
\[ z=z_1z_2+z_3\sqrt{z_1^2+1}\sqrt{z_2^2+1},\,\, z_3>1
\] 
учитывая из \cite[Гл.3, п.3.9.2]{Beitman_1} 
\begin{gather*}
\lim_{z_3 \to 1} \frac{e^{-i\,\pi\,\left( \mu -\frac{1}{2} \right)}Q_{n+v+\frac{1}{2}}^{\mu-\frac{1}{2}}(z_3)}{\sqrt{z_3^2-1}^{\mu -\frac{1}{2}}}=
\frac{\Gamma(n+v+\mu+1) \Gamma\left( \dfrac{1}{2}-\mu \right)}{2^{\mu+\frac{1}{2}}\Gamma( n+v-\mu+2)}
\\
\mu<\frac{1}{2}
\end{gather*}
а также вводя другие переменные вида
\begin{gather*}
z_1=\frac{r_1^2-1}{2r_1},\,\,z_2=\frac{1-r_2^2}{2r_2} \\
\text{где} \quad z=z_1z_2+\sqrt{z_1^2+1}\sqrt{z_2^2+1}=\frac{r_1^2+r_2^2}{2r_1r_2}
\end{gather*}
и с учетом того, что $ C_n^{v+1}(-x)=(-1)^nC_n^{v+1}(x) $, получим
\begin{multline} \label{Q_CC}
\frac{e^{-i\,\pi\,\mu}Q_{v}^{\mu}(z)}{\sqrt{z^2-1}^{\mu}}=
\frac{2^{4v+\mu+3} \Gamma\left( \dfrac{1}{2}-\mu \right) \Gamma(v+1)^2 (r_1r_2)^{v+\mu+1}}{\sqrt{\pi}\left( (r_1^2+1) (r_2^2+1) \right)^{v+\mu+1}}
\times
\\
\times
\sum_{n=0}^{\infty}\frac{ n! \, (n+v+1)\Gamma(n+v+\mu+1)}{\Gamma(n+v-\mu+2)\Gamma(n+2v+2)}
C_n^{v+1}\left( \frac{r_1^2-1}{r_1^2+1} \right)C_n^{v+1}\left( \frac{r_2^2-1}{r_2^2+1} \right)
\end{multline}
\[ \mu<\frac{1}{2} \]
или при $ v=l+\dfrac{M-3}{2},\,\, \mu=\dfrac{\nu-M+1}{2} $ из \eqref{eq:e13}
\begin{multline*}
\frac{e^{-i\,\pi\,\left( \frac{\nu-M+1}{2} \right)}Q_{l+\frac{M-3}{2}}^{\frac{\nu-M+1}{2}}(z)}{\sqrt{z^2-1}^{\frac{\nu-M+1}{2}}}=
\frac{2^{4l+\frac{3M+\nu-5}{2}} \Gamma\left( \dfrac{M-\nu}{2} \right) \Gamma\left( l+\dfrac{M-1}{2} \right)^2 (r_1r_2)^{l+\tfrac{\nu}{2}}}{\sqrt{\pi}\left( (r_1^2+1) (r_2^2+1) \right)^{l+\tfrac{\nu}{2}}}
\times
\\
\times
\sum_{n=0}^{\infty}\frac{ n! \, \left( n+l+\dfrac{M-1}{2} \right)\Gamma\left( n+l+\dfrac{\nu}{2} \right)}{\Gamma\left(n+l+M-\dfrac{\nu}{2} \right)\Gamma(n+2l+M-1)}
\times
\\
\times
C_n^{l+\tfrac{M-1}{2}}\left( \frac{r_1^2-1}{r_1^2+1} \right)
C_n^{l+\tfrac{M-1}{2}}\left( \frac{r_2^2-1}{r_2^2+1} \right)
\end{multline*}
\[ \text{при} \quad \nu<M \]
Отметим, что выражение \eqref{Q_CC} при $ \mu=0 $ эквивалентна выражению \cite[(11)]{Ossicini52} и совпадает с  \cite[(3.30)]{Radoslaw_1111.1661} (при $ r_s\to \dfrac{p_s}{q} $ ввиду однородности). Так как многочлены Гегенбауэра образуют ортогональную систему \cite[Гл.3, п.3.15(16)б(17)]{Beitman_1}, введем функции вида \eqref{eq:e_eta}, и тогда
\begin{multline*}
\frac{e^{-i\,\pi\,\left( \frac{\nu-M+1}{2} \right)}Q_{l+\frac{M-3}{2}}^{\frac{\nu-M+1}{2}}(z)}{\sqrt{z^2-1}^{\frac{\nu-M+1}{2}}}=
2^{\frac{\nu-M-1}{2}}\,\sqrt{\pi}\,\Gamma\left( \frac{M-\nu}{2} \right) (r_1r_2)^{\frac{\nu}{2}}
\left( (r_1^2+1) (r_2^2+1) \right)^{\frac{M-\nu}{2}}
\times
\\
\times
\sum_{n=0}^{\infty} \frac{\Gamma\left(n+l+\dfrac{\nu}{2} \right)}{\Gamma\left(n+l+M-\dfrac{\nu}{2} \right)}
\eta_{n,l}^{(M)}(r_1)\eta_{n,l}^{(M)}(r_2)
\end{multline*}
Таким образом \eqref{eq:e13} запишем как
\begin{multline*}
\int\limits_{0}^{\infty}du \,u^{\nu-1}
\frac{J_{l+\frac{M}{2}-1}(2r_1u)}{(r_1u)^{\frac{M}{2}-1}}\frac{J_{l+\frac{M}{2}-1}(2r_2u)}{(r_2u)^{\frac{M}{2}-1}}=
\\
=\frac{1}{4}\,
\left( (r_1^2+1) (r_2^2+1) \right)^{\frac{M-\nu}{2}}
\sum_{n=0}^{\infty} \frac{\Gamma\left(n+l+\dfrac{\nu}{2} \right)}{\Gamma\left(n+l+M-\dfrac{\nu}{2} \right)}
\eta_{n,l}^{(M)}(r_1)\eta_{n,l}^{(M)}(r_2)
\end{multline*}
учитывая это выражение получим для \eqref{eq:e12} в виде
\begin{multline*}
\int\limits_{0}^{\infty}\!\!du \,u^{\nu-1} 
\xi_{n_1,l}^{(\nu,M)}(u)\xi_{n_2,l}^{(\nu,M)}(u)
=(-1)^l\,\pi^{M} \sum_{n=0}^{\infty} \frac{\Gamma\left(n+l+\dfrac{\nu}{2} \right)}{\Gamma\left(n+l+M-\dfrac{\nu}{2} \right)}
\times
\\
\times
\int_E dr_1\,\rho(r_1)\phi_{n_1}(r_1)(r_1^2+1)^{\tfrac{M-\nu}{2}}\eta_{n,l}^{(M)}(r_1) 
\int_E dr_2\,\rho(r_2)\phi_{n_2}(r_2)(r_2^2+1)^{\tfrac{M-\nu}{2}}\eta_{n,l}^{(M)}(r_2)
\end{multline*}
Очевидно, что для
\begin{equation} \label{eq:e11}
\phi_n(r)=(r^2+1)^{\tfrac{M-\nu}{2}}\,\eta_{n,l}^{(M)}(r), \quad \rho(r)=\frac{r^{M-1}}{\left( r^2+1 \right)^{M-\nu}}
\end{equation}
при условии \eqref{eq:e_int_eta} и области $ E $ на всей действительной положительной оси, функция $ \xi_{n,l}^{(\nu,M)}(u) $ будет образовывать ортогональную систему \eqref{eq:c28}. 
Сопоставляя \eqref{1_bR} с \eqref{1_b11} а также с \eqref{eq:e8}-\eqref{eq:e10}, вводя функции \eqref{eq:c19} и \eqref{eq:eH}, получим \eqref{eq:c20}.

\subsection{ Функция $ \xi_{n,l}^{(\nu,M)}(u) $ и другие его представления.} \label{sec_1_3_2} 
Из \eqref{eq:e9} и \eqref{eq:e11} получим, что
\begin{equation} \label{eq:d6}
\xi_{n,l}^{(\nu,M)}(u)=
2 \, i^l \, \pi^{\frac{M}{2}}
\int\limits_{0}^{\infty}\!\!dr\,r^{M-1}
\frac{\eta_{n,l}^{(M)}(r)}{\left(r^2+1 \right)^{\frac{M-\nu}{2}}}
\frac{J_{l+\frac{M}{2}-1}(2ru)}{(ru)^{\frac{M}{2}-1}}
\end{equation}
или в виде, используя выражение для $ \eta_{n,l}^{(M)}(r) $ из \eqref{eq:e_eta}
\begin{multline} \label{eq:e20}
\xi_{n,l}^{(\nu,M)}(u)=
\frac{i^l\, 4 \, \pi^{\tfrac{M}{2}}}{\Gamma\left( l+\dfrac{M}{2} \right)}
\sqrt{\frac{\left(n+l+\displaystyle \frac{M-1}{2}\right)\Gamma\left(n+2l+M-1\right)}{ \Gamma\left(n+1\right)}}
\times
\\
\times
u^{l+M-\nu}\int\limits_{0}^{\infty}dt \,
\frac{t^{l+\frac{M}{2}}J_{l+\frac{M}{2}-1}(2t)}{(t^2+u^2)^{l+M-\tfrac{\nu}{2}}}
\,{}_{2}F_{1} \left[ 
\left.
\begin{gathered}
{ -n \quad n+2l+M-1 } \\
{ l+\frac{M}{2}  }
\end{gathered}
\right| \frac{u^2}{t^2+u^2} \right]
\end{multline}
В отличии от $ \eta_{n,l}^{(M)}(r) $  \eqref{eq:e_eta}, функция $ \xi_{n,l}^{(\nu,M)}(u) $ не выражается элементарно в общем случае. Так, раскрывая гипергеометрический ряд в \eqref{eq:e20} и используя интегралы от функции Бесселя \cite[Гл.2, п.2.12.4(28)]{Prudnikov_2} 
\begin{equation*}
\int\limits_{0}^{\infty}\!\!dx\, \frac{x^{v+1}J_v(cx)}{(x^2+z^2)^{\rho}}=
\frac{c^{\rho-1}z^{v-\rho+1}}{2^{\rho-1}\Gamma( \rho )}K_{v-\rho+1}(cz), \quad
-1<v<2\rho-\frac{1}{2}
\end{equation*}
представим \eqref{eq:e20} в виде ряда (здесь по определению для функции Макдональдса выполняется условие $ K_{-\mu}(z)=K_{\mu}(z) $)

\begin{multline*}
\xi_{n,l}^{(\nu,M)}(u)
=\frac{i^l\, 4 \, \pi^{\tfrac{M}{2}}}{\Gamma\left( l+\dfrac{M}{2} \right) \Gamma\left(l+M-\displaystyle\frac{\nu}{2}\right)}
\sqrt{\frac{\left(n+l+\displaystyle \frac{M-1}{2}\right)\Gamma\left(n+2l+M-1\right)}{\Gamma\left(n+1\right)}} 
\times
\\
\times
\sum_{m=0}^{n}
\frac{(-n)_{m} (n+2l+M-1)_m}{\left( l +\displaystyle \frac{M}{2} \right)_m \left( l + M - \displaystyle \frac{\nu}{2} \right)_m m!} 
(u)^{l+m+\frac{M-\nu}{2}} K_{\frac{M-\nu}{2}+m}(2u)
\end{multline*}
из \cite[Гл.2, п.2.3.16(1)]{Prudnikov_1} ($ u^2>0 $)
\begin{equation*}
\int\limits_{0}^{\infty}\!\!dt \,e^{-t-\frac{u^2}{t}} t^{\frac{M-\nu}{2}+m-1}=
2(u)^{\frac{M-\nu}{2}+m} K_{\frac{M-\nu}{2}+m}(2u)
\end{equation*}
Сопоставляя это выражение с предыдущем, получим\eqref{eq:c16}. Аналогично выражению \\ \eqref{eq:e20} можно также получить общее выражение

\begin{multline*}
\xi_{n,l}^{(\nu,M)}(u)=
\frac{i^l\, 4 \, \pi^{\tfrac{M}{2}}\,\Gamma\left( \dfrac{M{-}\nu}{2}{+}\lambda{+}1 \right)}{\Gamma\left( l{+}\dfrac{M}{2} \right) \Gamma\left(l{+}M{-}\dfrac{\nu}{2}\right)}
\sqrt{\frac{\left(n{+}l{+}\dfrac{M{-}1}{2}\right)\Gamma\left(n{+}2l{+}M{-}1\right)}{ \Gamma\left(n{+}1\right)}}
\times
\\
\times
u^{l+M-\nu}\int\limits_{0}^{\infty}\!\!dt \,
\frac{t^{\lambda+1}J_{\lambda}(2t)}{(t^2{+}u^2)^{\tfrac{M-\nu}{2}+\lambda+1}}
\,{}_{3}F_{2} \left[ 
\left.
\begin{gathered}
{ -n \quad n{+}2l{+}M{-}1 \quad \frac{M{-}\nu}{2}{+}\lambda{+}1 } \\
{ l{+}\frac{M}{2} \quad l{+}M{-}\frac{\nu}{2} }
\end{gathered}
\right| \frac{u^2}{t^2{+}u^2} \right]
\end{multline*}
\[ \lambda>-1,\,\, \lambda+M-\nu+\frac{3}{2}>0 \]

\begin{center}
	\textbf{Примеры}
\end{center}
При $ N=2 $ для $ M $--мерных векторов $ \mathbf{r}_s $ из \eqref{eq:c20} с учетом \eqref{eq:c28} и \eqref{int_Y} получим
\begin{multline*}
\frac{1}{\left| \mathbf{r}_1 - \mathbf{r}_2 \right|^{\nu}}=
\frac{\pi^{\tfrac{M}{2}}\Gamma\left(\dfrac{M-\nu}{2}\right)}{\Gamma\left( \dfrac{\nu}{2}\right)} 
\left(\left(r_1^2+1\right) \left(r_2^2+1\right) \right)^{\tfrac{M-\nu}{2}} 
\times 
\\
\times 
\sum_{n=0}^{\infty} \sum_{l=0}^{\infty} \sum_{\mathbf{m}_k}
\frac{\Gamma\left(l+n+\dfrac{\nu}{2} \right)}{\Gamma\left(l+n+M-\dfrac{\nu}{2}\right)} 
H_{n,l,\mathbf{m}_k}^{(M)}(\mathbf{r}_1) 
H_{n,l,\mathbf{m}_k}^{(M)}(\mathbf{r}_2)^{\ast}
\end{multline*}
а также
\begin{multline}
\label{1_b015}
\frac{1}{\left| \mathbf{r}_1 + \mathbf{r}_2 \right|^{\nu}}=
\frac{\pi^{\tfrac{M}{2}}\Gamma\left(\dfrac{M-\nu}{2}\right)}{\Gamma\left( \dfrac{\nu}{2}\right)} 
\left(\left(r_1^2+1\right) \left(r_2^2+1\right) \right)^{\tfrac{M-\nu}{2}} 
\times 
\\
\times 
\sum_{n=0}^{\infty} \sum_{l=0}^{\infty} \sum_{\mathbf{m}_k}
\frac{(-1)^l\,\Gamma\left(l+n+\dfrac{\nu}{2} \right)}{\Gamma\left(l+n+M-\dfrac{\nu}{2}\right)} 
H_{n,l,\mathbf{m}_k}^{(M)}(\mathbf{r}_1) 
H_{n,l,\mathbf{m}_k}^{(M)}(\mathbf{r}_2)^{\ast}
\end{multline}
Как было сказано во введении \ref{I}, для $ M $-мерных векторов $ \mathbf{r} $ и для соответствующих их гиперсферических функции $ Y_{l,\mathbf{m_k}}(\boldsymbol{\zeta}) $ существуют $ \dfrac{(2M-2)!}{(M-1)!\,M!} $  эквивалентных представлении. Соответственно столько же эквивалентных представлении может быть и для этих примеров и для \eqref{eq:c20}     

\section{Связь между функциями \eqref{eq:c19} и \eqref{eq:eH}} 
\noindent \textbf{4.1} Аналогично выражению \eqref{eq:d6}, где установлена связь между $ \xi $ и $ \eta $ функциями, можно показать что
\begin{equation}\label{4_1}
\Xi_{n,l,\mathbf{m_k}}(u,\boldsymbol{\zeta}_u)=
\int\limits\!\! d\mathbf{r} \,\frac{H_{n,l,\mathbf{m_k}}^{(M)}(\mathbf{r})^{\ast}}{(r^2+1)^{\frac{M-\nu}{2}}}e^{i\,2(\mathbf{r}\mathbf{u})}
\end{equation} 
где $ \mathbf{u}=u\boldsymbol{\zeta}_u $. Действительно, из \cite[Гл.10, п.10.20(7)]{Beitman_2}
\begin{equation*}
e^{ixy}=\Gamma( \lambda )\left( \frac{2}{y} \right)^{\lambda}\sum_{l=0}^{\infty} i^l(l+\lambda) J_{l+\lambda}(y) C_l^{\lambda}(x), \quad -1<x<1,\,\, \lambda>0
\end{equation*}
при $ (\mathbf{r}\mathbf{u})=ru(\boldsymbol{\zeta}_r\boldsymbol{\zeta}_u) $ если взять $ y=ru, \,\,x=(\boldsymbol{\zeta}_r\boldsymbol{\zeta}_u),\,\,\lambda=\dfrac{M}{2}-1 $ и используя теорему сложения для многочленов Гегенбауэра \eqref{th_C} получим разложение для плоской волны в $ M $--мерном пространстве 
\begin{equation*}
e^{i(\mathbf{r}\mathbf{u})}=(2\pi)^{\tfrac{M}{2}}\sum_{l,\mathbf{m_k}} i^l\,\frac{J_{l+\tfrac{M}{2}-1}(ru)}{(ru)^{\tfrac{M}{2}-1}} Y_{l,\mathbf{m_k}}(\boldsymbol{\zeta}_r) Y_{l,\mathbf{m_k}}(\boldsymbol{\zeta}_u)^{\ast}
\end{equation*}
подставляя это в выражение в \eqref{4_1} 
\begin{multline*}
\int\limits\!\! d\mathbf{r} \,\frac{H_{n,l,\mathbf{m_k}}^{(M)}(\mathbf{r})^{\ast}}{(r^2+1)^{\frac{M-\nu}{2}}} e^{i\,2(\mathbf{r}\mathbf{u})}=
(2\pi)^{\tfrac{M}{2}}\int\limits_{0}^{\infty}\!\! r^{M-1} dr \!\!\int\limits\!\! d\Omega_{\boldsymbol{\zeta}_r} \,
\frac{\eta_{n,l}^{(M)}(r)Y_{l,\mathbf{m_k}}(\boldsymbol{\zeta}_r)^{\ast}}{(r^2+1)^{\frac{M-\nu}{2}}}
\times
\\
\times
\sum_{l',\mathbf{m_k}'} i^{l'}\,\frac{J_{l'+\tfrac{M}{2}-1}(2ru)}{(2ru)^{\tfrac{M}{2}-1}} Y_{l',\mathbf{m_k}'}(\boldsymbol{\zeta}_r) Y_{l',\mathbf{m_k}'}(\boldsymbol{\zeta}_u)^{\ast}
\end{multline*}
интегрируя по $ d\Omega_{\boldsymbol{\zeta}_r} $ по всему $ M-1 $--мерному пространству, с учетом \eqref{int_Y} и снимая суммы по $ l',\mathbf{m_k}' $ получим для последнего выражения
\begin{equation*}
i^l\,(2\pi)^{\tfrac{M}{2}}\int\limits_{0}^{\infty}\!\!dr\, r^{M-1}  \,
\frac{\eta_{n,l}^{(M)}(r)}{(r^2+1)^{\frac{M-\nu}{2}}}
\frac{J_{l+\tfrac{M}{2}-1}(2ru)}{(2ru)^{\tfrac{M}{2}-1}} Y_{l,\mathbf{m_k}}(\boldsymbol{\zeta}_u)^{\ast} 
\end{equation*}
С учетом \eqref{eq:d6} и \eqref{eq:c19} мы получим \eqref{4_1}.

\medskip
\noindent \textbf{4.2} Покажем, что вид выражения \eqref{eq:c20} не меняется при последовательных пределах \\  
$ \mathbf{r}_{_N} \to 0, \, \mathbf{r}_{_{N-1}} \to 0, \, \ldots,\mathbf{r}_{_{3}} \to 0  $,
и при $ N=2 $  она будет соответствовать виду \eqref{1_b015}. Так,
\begin{multline}
\label{eq:c21}
\frac{1}{\left|\mathbf{r}_1+\ldots+\mathbf{r}_{_{N-1}}\right|^{\nu}}=
\lim\limits_{\mathbf{r}_{_N} \to 0} \frac{1}{\left|\mathbf{r}_1+\ldots+\mathbf{r}_{_N}\right|^{\nu}}= \\
=\frac{2\Gamma\left(\dfrac{M-\nu}{2}\right) }
{S_{_M} \, \Gamma\left(\dfrac{M}{2}\right)\Gamma\left(\dfrac{\nu}{2}\right)}
\sum_{\mathbf{n',l',m_k'}} \,\, \prod_{p=1}^{N-1}
\left( \left(r_p^2+1\right)^{\tfrac{M-\nu}{2}} H_{n_p,l_p,\mathbf{m}_{\mathbf{k}}^{(p)}}^{(M)}(\mathbf{r}_p)\right)
\times
\\
\times
\int\limits_{0}^{\infty}du \, u^{\nu-1} \int d\Omega_{\boldsymbol{\zeta}} \prod_{p=1}^{N-1}
\Xi_{n_p,l_p,\mathbf{m}_{\mathbf{k}}^{(p)}}^{(\nu,M)}(u,\boldsymbol{\zeta}) 
\times
\\
\times
\lim\limits_{\mathbf{r}_{_N} \to 0} \, \sum_{n_{_N}=0}^{\infty} \, \sum_{l_{_N}=0}^{\infty} \, \sum_{\mathbf{m}_{\mathbf{k}}^{(N)}}
\left(r_{_N}^2+1\right)^{\tfrac{M-\nu}{2}} H_{n_{_N},l_{_N},\mathbf{m}_{\mathbf{k}}^{(N)}}^{(M)}(\mathbf{r}_{_N}) \,
\Xi_{n_{_N},l_{_N},\mathbf{m}_{\mathbf{k}}^{(N)}}^{(\nu,M)}(u,\boldsymbol{\zeta})       
\end{multline}
где $ \sum\limits_{\mathbf{n',l',m_k'}} $  означает все суммы кроме по $ n_{_N},l_{_N},\mathbf{m}_{\mathbf{k}}^{(N)} $ . Такое возможно, если при любых $ (u,\boldsymbol{\zeta}) $ , выполняется условие.
 \begin{equation}
	\label{eq:c22}
	\lim\limits_{\mathbf{r}_{_N} \to 0} \, \sum_{n_{_N}=0}^{\infty} \, \sum_{l_{_N}=0}^{\infty} \, \sum_{\mathbf{m}_{\mathbf{k}}^{(N)}}
	\left(r_{_N}^2+1\right)^{\tfrac{M-\nu}{2}} H_{n_{_N},l_{_N},\mathbf{m}_{\mathbf{k}}^{(N)}}^{(M)}(\mathbf{r}_{_N}) \,
	\Xi_{n_{_N},l_{_N},\mathbf{m}_{\mathbf{k}}^{(N)}}^{(\nu,M)}(u,\boldsymbol{\zeta}) =1
\end{equation}
Сумму по $ \mathbf{m}_{\mathbf{k}} $  снимем с учетом теоремы сложения \eqref{th_C} и из определения \eqref{eq:c19} и \eqref{eq:eH} получим
\begin{multline}
\label{eq:c22_1}
\frac{1}{S_{_M}}\lim\limits_{r_{_N} \to 0} \left(r_{_N}^2+1\right)^{\tfrac{M-\nu}{2}} \sum_{n_{_N}=0}^{\infty} \, \sum_{l_{_N}=0}^{\infty} 
\eta_{n_{_N},l_{_N}}^{(M)}(r_{_N}) \,
\xi_{n_{_N},l_{_N}}^{(\nu,M)}(u) 
\times
\\
\times
\frac{l_{_N}+\dfrac{M}{2}-1}{\dfrac{M}{2}-1}
C_{l_{_N}}^{\tfrac{M}{2}-1} \left( ( \boldsymbol{\zeta}^{(N)} \boldsymbol{\zeta} ) \right) 
=1
\end{multline}
Из \eqref{eq:e_eta} легко получить, что
\begin{equation} \label{eq:a5}
\eta_{n,l}^{\left(M\right)}\left(0\right)=\delta_{l,0}\frac{2\left(-1\right)^n}{\Gamma\left(\dfrac{M}{2}\right)}
\sqrt{\frac{\left(n+\dfrac{M-1}{2}\right)\Gamma\left(n+M-1\right)}{n!}} 
\end{equation} 
таким образом в \eqref{eq:c22_1} при $ r_{_N} \to 0 $  в суммах по $ l_{_N} $ отличны от нуля будут только члены при $ l_{_N}=0 $ , и в итоге 
\begin{equation} \label{eq:c22_2}
\frac{1}{S_{_M}}
\sum_{n_{_N}=0}^{\infty} \,
\eta_{n_{_N},0}^{(M)}(0) \,
\xi_{n_{_N},0}^{(\nu,M)}(u) =1
\end{equation}
Это выражение можно получить из условия ортогональности \eqref{eq:c28}. Или также аналогично с учетом \eqref{int_Y} условием ортогональности функции \eqref{eq:c19}
\begin{multline}
\label{eq:c27_1}
\int\limits_{0}^{\infty}\!\!du \, u^{\nu-1} \int\!\! d\Omega_{\boldsymbol{\zeta}} \,
\Xi_{n_1,l_1,\mathbf{m}_{\mathbf{k}}^{(1)}}^{(\nu,M)}(u,\boldsymbol{\zeta})
\Xi_{n_2,l_2,\mathbf{m}_{\mathbf{k}}^{(2)}}^{(\nu,M)}(u,\boldsymbol{\zeta})^{\ast} =
\\
=\frac{\pi^{M} \Gamma\left(\dfrac{\nu}{2}+l_2+n_2\right) }
{\Gamma\left(M-\dfrac{\nu}{2}+l_2+n_2\right) }
\delta_{n_1,n_2} \delta_{l_1,l_2} \delta_{\mathbf{m}_{\mathbf{k}}^{(1)},\mathbf{m}_{\mathbf{k}}^{(2)}}       \end{multline} 
можно получить выражение \eqref{eq:c22} если учесть что при любых $ \boldsymbol{\zeta}, \,\, Y_{0,0}(\boldsymbol{\zeta})\sqrt{S_{_M}}=1 $ .  \\ В самом деле выражение \eqref{eq:c22_2} есть разложение постоянной единица по ортогональным функциям \eqref{eq:c16}, с коэффициентами $ \eta_{n_{_N},0}^{(M)}(0) / S_{_M} $ . Действительно, беря в левой и правой части выражения \eqref{eq:c22_2} интеграл по $ \xi_{n,0}^{(\nu,M)}(u) $  с весом $ u^{\nu-1} $, получим 
\begin{equation*}
\frac1{S_{_M}} \int\limits_{0}^{\infty}du \, u^{\nu-1} \xi_{n,0}^{(\nu,M)}(u)
\sum_{n_{_N}=0}^{\infty} \eta_{n_{_N},0}^{(M)}(0) \xi_{n_{_N},0}^{(\nu,M)}(u)
=\int\limits_{0}^{\infty}du \, u^{\nu-1} \xi_{n,0}^{(\nu,M)}(u)
\end{equation*}
и учитывая \eqref{eq:c28} получим
\begin{multline*}
\frac1{S_{_M}} \sum_{n_{_N}=0}^{\infty} \eta_{n_{_N},0}^{(M)}(0)
\frac{\pi^M \Gamma\left(\dfrac{\nu}{2}+n_{_N}\right) }{\Gamma\left(M-\dfrac{\nu}{2}+n_{_N}\right)}
\delta_{n_{_N},n}=
\eta_{n,0}^{(M)}(0)
\frac{\pi^M \Gamma\left(\dfrac{\nu}{2}+n\right) }{S_{_M}\Gamma\left(M-\dfrac{\nu}{2}+n\right)}=
\\
=\int\limits_{0}^{\infty}du \, u^{\nu-1} \xi_{n,0}^{(\nu,M)}(u)
\end{multline*}
Используя значение $ \eta_{n,0}^{(M)}(0) $ из  \eqref{eq:a5} запишем последнее выражение в виде,
\begin{equation*}
\int\limits_{0}^{\infty}du \, u^{\nu-1} \xi_{n,0}^{(\nu,M)}(u)
=\frac{(-1)^n \pi^{\tfrac{M}{2}}\Gamma\left(\dfrac{\nu}{2}+n\right) }
{\Gamma\left(M-\dfrac{\nu}{2}+n\right)}
\sqrt{\frac{\left(n+\dfrac{M-1}{2}\right)\Gamma(n+M-1)}{n!}}
\end{equation*}
С другой стороны из \eqref{eq:c16} с учетом того, что
\begin{equation*}
\int\limits_{0}^{\infty}\!\! du \,u^{\nu-1}e^{-\frac{u^2}{z}}=
\frac{1}{2}z^{\tfrac{\nu}{2}}\,\Gamma\left( \frac{\nu}{2} \right)
\end{equation*} 
получим
\begin{multline*}
	\int\limits_{0}^{\infty}du \, u^{\nu-1} \xi_{n,0}^{(\nu,M)}(u)=
	\frac{\pi^{\tfrac{M}{2}}\Gamma\left(\dfrac{\nu}{2}\right) }
	{\Gamma\left(\dfrac{M}{2}\right) \Gamma\left(M-\dfrac{\nu}{2}\right)}
	\sqrt{\frac{\left(n+\dfrac{M-1}{2}\right)\Gamma(n+M-1)}{n!}}
	\times
	\\
	\shoveright{
		\times
		\int\limits_{0}^{\infty}\!\!dz \, e^{-z} z^{\tfrac{M}{2}-1}
		\,{}_2F_2 \left[ \left.
		\begin{gathered}
		{ -n \quad n+M-1  } \\
		{ \frac{M}{2} \quad M-\dfrac{\nu}{2} }
		\end{gathered}
		\right| z   \right]=}
	\\
	=\frac{\pi^{\tfrac{M}{2}}\Gamma\left(\dfrac{\nu}{2}\right) }
	{ \Gamma\left(M-\dfrac{\nu}{2}\right)}
	\sqrt{\frac{\left(n+\dfrac{M-1}{2}\right)\Gamma(n+M-1)}{n!}}
	\,{}_2F_1 \left[ \left.
	\begin{gathered}
	{ -n \quad n+M-1  } \\
	{ M-\dfrac{\nu}{2} }
	\end{gathered}
	\right| 1   \right]=
	\\
	=\frac{\pi^{\tfrac{M}{2}}\Gamma\left(\dfrac{\nu}{2}\right) }
	{ \Gamma\left(M-\dfrac{\nu}{2}\right)}
	\sqrt{\frac{\left(n+\dfrac{M-1}{2}\right)\Gamma(n+M-1)}{n!}}
	\frac{\left(-\dfrac{\nu}{2}-n+1\right)_n }{\left(M-\dfrac{\nu}{2}\right)_n}= 
	\\
	=\frac{(-1)^n \pi^{\tfrac{M}{2}}\Gamma\left(\dfrac{\nu}{2}+n\right) }
	{\Gamma\left(M-\dfrac{\nu}{2}+n\right)}
	\sqrt{\frac{\left(n+\dfrac{M-1}{2}\right)\Gamma(n+M-1)}{n!}}
\end{multline*}
что и доказывает \eqref{eq:c22_2} а следовательно и \eqref{eq:c22} , \eqref{eq:c22_1}. 

Аналогично обратному преобразованию \eqref{eq:c21} можно получить соотношение вида \\ (теорему сложения)
\begin{multline} \label{eq:c33}
\left(\left| \mathbf{r}_1+\ldots+\mathbf{r}_{_N}\right|^2+1  \right)^{\tfrac{M-\nu}{2}}
H_{n',l',\mathbf{m_k}'}^{(M)}(\mathbf{r}_1+\ldots+\mathbf{r}_{_N})= 
\\
=
\frac{\Gamma\left(M-\dfrac{\nu}{2}+l+n\right) }
{\pi^M \Gamma\left(\dfrac{\nu}{2}+l+n\right)}
\sum_{\mathbf{n,l},\mathbf{m}_{\mathbf{k}}^{(i)}} \, \prod_{p=1}^{N}\left(r_p^2+1\right)^{\tfrac{M-\nu}{2}} 
H_{n_p,l_p,\mathbf{m}_{\mathbf{k}}^{(p)}}^{(M)}(\mathbf{r}_p) 
\times
\\
\times  
\int\limits\limits_{0}^{\infty}\!\!du \!\! \int\!\! d\Omega_{\boldsymbol{\zeta}} \, u^{\nu-1}  \,
\Xi_{n',l',\mathbf{m_k}'}^{(\nu,M)}(u,\boldsymbol{\zeta}) \prod_{s=1}^{N} 
\Xi_{n_s,l_s,\mathbf{m}_{\mathbf{k}}^{(s)}}^{(\nu,M)}(u,\boldsymbol{\zeta})
\end{multline}
которая легко выводится из \eqref{eq:c20}   при
\begin{equation*}
\lim\limits_{\mathbf{y} \to \mathbf{r}_{k+1}+\ldots+\mathbf{r}_{_N}} \quad
\frac{1}{\left| \mathbf{r}_1+\ldots+\mathbf{r}_{k}+\mathbf{y}\right|^{\nu}}=
\frac{1}{\left| \mathbf{r}_1+\ldots+\mathbf{r}_{_N}\right|^{\nu}}
\end{equation*}
и используя \eqref{eq:c22} и \eqref{eq:c27_1}. 

Отметим, что в \eqref{eq:c33} суммирование производится по всем $ n_q,l_q,\mathbf{m_k}^{(q)}, \,\,q=1\ldots N $, кроме $ n',l',\mathbf{m_k}' $. Выражение \eqref{eq:c33} имеет особый интерес в физических задачах и приложениях.

\appendix 
\section{Приложение} \label{append}
В этом приложении дан вывод разложения \eqref{1_f6}. По сути здесь представлен аналогичная часть нашей работы \cite{AkhR_UNC}, и поэтому этот раздел мы обозначили как приложение.  
Так, исходя из соотношения \cite[Гл.6, п.6.7.3(3)]{Prudnikov_3} 
\begin{equation} \label{1_v1}
\frac{1}{\left(h^2-2h\omega+1\right)^a}=\sum_{l=0}^{\infty} \frac{\left(a\right)_l}{\left(b-1\right)_l} h^l C_l^{b-1} \left(\omega\right)  {}_2F_1 \left[ \left.\begin{matrix}
{ a+l \quad a-b+1 } \\
{ b+l }
\end{matrix} \right| h^2 \right], \quad (a>-1)
\end{equation}
если положить $ a{=}\dfrac{\nu}{2},\,\omega{=}\cos\omega_{12}{=}\dfrac{(\mathbf{r}_1\mathbf{r}_2)}{r_1r_2} $ --- косинус угла между векторами, и $ h{=}\dfrac{r_1}{r_2}<1 $, можно получить из \eqref{1_v1}   
\begin{multline} \label{1_a1}
\frac{1}{\left| \mathbf{r}_1 - \mathbf{r}_2 \right|^{\nu}}= 
\frac{1}
{r_2^{\nu}} 
\sum_{l=0}^{\infty} 
\frac{\left(\dfrac{\nu}{2}\right)_l}{\left(b-1\right)_l} 
\left( \frac{r_1}{r_2} \right)^l 
C_l^{b-1} \left( \cos \omega_{12} \right) 
{}_2F_1 \left[ \left.\begin{matrix}
{  l+\dfrac{\nu}{2}  \quad 1-b+\dfrac{\nu}{2}  } \\
{  b+l  }
\end{matrix} \right| 
\left( \frac{r_1}{r_2} \right)^2 \right]  = 
\\
= \frac{1}{r_2^{\nu}} 
\sum_{l=0}^{\infty} 
\frac{C_l^{b-1} \left( \cos \omega_{12} \right)}{\left(b-1\right)_l}
\sum_{\mu =0}^{\infty} \frac{\left(\dfrac{\nu}{2}\right)_{l+\mu} \left( 1-b+\dfrac{\nu}{2}\right)_{\mu}}{\left(b+l\right)_{\mu} \mu!}
\left( \frac{r_1}{r_2} \right)^{l+2\mu}
\end{multline}
Используя представление функции Лежандра второго рода через гипергеометрические функции Гаусса из \cite[Гл.3, п.3.2(45)]{Beitman_1}
\begin{equation*}
e^{-i\pi\mu}Q_v^{\mu}(z)=
\frac{\sqrt{\pi}\,2^{\mu}(z^2-1)^{\frac{\mu}{2}}}{\left( z+\sqrt{z^2-1} \right)^{v+\mu+1}}
\frac{\Gamma( v+\mu+1 )}{\Gamma\left( v+\dfrac{3}{2} \right)}
\,{}_{2}F_{1} \left[ 
\left.
\begin{gathered}
{ \mu+\frac{1}{2} \quad v+\mu+1 } \\
{ v+\frac{3}{2} }
\end{gathered}
\right|  \frac{z-\sqrt{z^2-1}}{z+\sqrt{z^2-1}}  \right]
\end{equation*}
получим при 
\begin{gather} \label{1_v3}
v=l+b-\frac{3}{2}, \quad \mu=\frac{\nu+1}{2}-b \\
\label{1_v4}
z=\frac{r_1^2+r_2^2}{2\,r_1 r_2}>1
\end{gather}
для гипергеометрической функции Гаусса в \eqref{1_a1}
\begin{equation*}
\,{}_{2}F_{1} \left[ 
\left.
\begin{gathered}
{ l{+}\frac{\nu}{2} \quad 1{-}b{+}\frac{\nu}{2} } \\
{ b{+}l }
\end{gathered}
\right|  \left( \frac{r_1}{r_2} \right)^2  \right]=
\frac{\Gamma( b{+}l )}{\sqrt{\pi}\,2^{\tfrac{\nu+1}{2}-b}\,\Gamma\left( l{+}\dfrac{\nu}{2} \right)}
\left( \frac{r_2}{r_1} \right)^{l+\tfrac{\nu}{2}} \,\,
\frac{e^{-i\pi\left(\tfrac{\nu+1}{2}-b\right)}\,\, Q_{l+b-\tfrac32}^{\tfrac{\nu+1}{2}-b}(z) }
{\left(\sqrt{z^2{-}1}\right)^{\tfrac{\nu+1}{2}-b}}
\end{equation*}
и в итоге получим для \eqref{1_a1} выражение вида
\begin{multline} \label{1_v5}
\frac{1}
{\left| \mathbf{r}_1 - \mathbf{r}_2 \right|^{\nu}}=
\frac{1}
{(r_1 r_2)^{\tfrac{\nu}{2}} \, \sqrt{\pi}\, \Gamma\left(\dfrac{\nu}{2}\right) 2^{\tfrac{\nu +1}{2}-b}}
\times
\\
\times     
\sum_{l=0}^{\infty}
(b+l-1)\Gamma(b-1) C_l^{b-1} \left( \cos(\omega_{12}) \right) 
\frac{e^{-i\pi\left(\tfrac{\nu+1}{2}-b\right)} Q_{l+b-\tfrac32}^{\tfrac{\nu+1}{2}-b}(z) }
{\left(\sqrt{z^2-1}\right)^{\tfrac{\nu+1}{2}-b}}
\end{multline}
Здесь и далее $ Q^{\mu}_{v}(z) $---функция Лежандра второго рода.

Детальное изложение  \eqref{1_a1} и \eqref{1_v5} и ее обсуждение в связи с другими подобными формулами и приложением для фундаментального  решения полигармонического уравнения Лапласа приведен в \cite{Cohl_2013_1105.2735}.  

Для того, чтобы найти такие функции $ f_n(r) $ и коэффициенты $ q_{n,m} $ при котором в \eqref{1_v5} можно представить разложение вида
\begin{equation*}
\frac{e^{-i\pi\mu}Q_{v}^{\mu}(z)}{\sqrt{z^2-1}^{\mu}}=
\sum_{n}q_{n,m}f_n(r_1)f_m(r_2)
\end{equation*}
(где $ z $ равно \eqref{1_v4}) ограничимся пока построением разложения при 
\begin{equation*}
z=z_1z_2+z_3\sqrt{z_1^2+1}\sqrt{z_2^2+1}
\end{equation*}  
по $ z_1,\,z_2\,z_3 $.
В соответствии с определением функции Лежандра второго рода представим его в виде ( $ z>1 $) \cite[Гл.3, п.3.2(5)]{Beitman_1}
\begin{multline} \label{1_f1}
\frac{e^{-i\pi\mu}Q_{v}^{\mu}(z)}{\sqrt{z^2-1}^{\mu}}=
\frac{\sqrt{\pi}\Gamma( v+\mu+1 )}{2^{v+1}z^{v+\mu+1}\Gamma\left( v+\dfrac{3}{2} \right)}
\,{}_{2}F_{1} \left[ 
\left.
\begin{gathered}
{ \frac{v+\mu+1}{2} \quad \frac{v+\mu+2}{2} } \\
{ v+\frac{3}{2} }
\end{gathered}
\right|  \frac{1}{z^2}  \right]=
\\
=\sqrt{\pi}\,2^{\mu}\sum_{n=0}^{\infty}\frac{\Gamma( v+\mu+1+2n )}{\Gamma\left( v+n+\dfrac{3}{2} \right)n!}\frac{1}{(2z)^{v+\mu+1+2n}}
\end{multline}
и положим в нем $ z=c_1+c_2 $ , а 
\[ \frac{1}{z^{\beta}}=\frac{1}{c_1^{\beta}\Gamma( \beta )}\sum_{m=0}^{\infty}(-1)^m\frac{\Gamma( \beta+m )}{m!}\left( \frac{c_1}{c_2} \right)^{\beta+m},\quad (\beta=v+\mu+1+2n) \] 
перепишем в виде суммы, которая сходится при любом $ \beta $  и $ \left| \dfrac{c_1}{c_2} \right| <1 $ .  Подставляя последнее выражение в \eqref{1_f1} и объединяя в сумме порядок по  $ \dfrac{c_1}{c_2} $,  запишем в виде
\begin{multline} \label{1_f2}
\frac{e^{-i\pi\mu}Q_{v}^{\mu}(c_1+c_2)}{\sqrt{(c_1+c_2)^2-1}^{\mu}}=
\frac{\sqrt{\pi}\,2^{\mu}}{(2c_2)^{v+\mu+1}\Gamma\left( v+\dfrac{3}{2} \right)}
\times
\\
\times
\sum_{n=0}^{\infty}\left( -\frac{c_1}{c_2} \right)^n \frac{\Gamma( v+\mu+1+n )}{n!}
\,{}_{2}F_{1} \left[ 
\left.
\begin{gathered}
{ -\frac{n}{2} \quad -\frac{n}{2}+\frac{1}{2} } \\
{ v+\frac{3}{2} }
\end{gathered}
\right|  \frac{1}{c_1^2}  \right]
\end{multline}
Отметим, что гипергеометрическая функция в этом выражении – конечный полином, поэтому никаких ограничений на $ c_1 $  не накладывается.
Далее представим гипергеометрическую функцию в \eqref{1_f2} в виде \cite[Гл.2, п.2.11(16)]{Beitman_1}
\begin{equation*}
{}_{2}F_{1} \left[ 
\left.
\begin{gathered}
{ -\frac{n}{2} \quad -\frac{n}{2}+\frac{1}{2} } \\
{ v+\frac{3}{2} }
\end{gathered}
\right|  \frac{1}{c_1^2}  \right]=
\left( \frac{\sqrt{c_1^2-1}}{c_1} \right)^n
\,{}_{2}F_{1} \left[ 
\left.
\begin{gathered}
{ -n \quad n+2v+2 } \\
{ v+\frac{3}{2} }
\end{gathered}
\right|  \frac{1}{2}\left( 1-\frac{c_1}{\sqrt{c_1^2-1}}\right)  \right]
\end{equation*}
и в свою очередь из \cite[Гл.3, п.3.15(3)]{Beitman_1}  
\begin{equation*}
{}_{2}F_{1} \left[ 
\left.
\begin{gathered}
{ -n \quad n+2v+2 } \\
{ v+\frac{3}{2} }
\end{gathered}
\right|  \frac{1}{2}\left( 1-\frac{c_1}{\sqrt{c_1^2-1} }\right)  \right]=
\frac{n!\Gamma( 2v+2 )}{\Gamma( n+2v+2 )}
C_n^{v+1}\left( \frac{c_1}{\sqrt{c_1^2-1}} \right)
\end{equation*}
Таким образом представим \eqref{1_f2} в виде
\begin{multline} \label{1_f3}
\frac{e^{-i\pi\mu}Q_{v}^{\mu}(c_1+c_2)}{\sqrt{(c_1+c_2)^2-1}^{\mu}}=
\frac{\sqrt{\pi}\,2^{\mu}\Gamma( 2v+2 )}{(2c_2)^{v+\mu+1}\Gamma\left( v+\dfrac{3}{2} \right)}
\times
\\
\times
\sum_{n=0}^{\infty}
\left( -\frac{\sqrt{c_1^2-1}}{c_2} \right)^n 
\frac{\Gamma( v+\mu+1+n )}{\Gamma( 2v+2+n )}
C_n^{v+1}\left( \frac{c_1}{\sqrt{c_1^2-1}} \right)
\end{multline}
Если положить
\begin{equation} \label{1_f4}
c_1=z_1z_2,\,\, c_2=z_3\sqrt{z_1^2+1}\sqrt{z_2^2+1},\quad z_3>1
\end{equation}
и ввести  вспомогательные переменные
\begin{equation} \label{1_f5} 
y_1=\frac{z_1}{\sqrt{z_1^2+1}},\quad y_2=\frac{z_2}{\sqrt{z_2^2+1}}
\end{equation}
то
\begin{equation*}
\left( -\frac{\sqrt{c_1^2-1}}{c_2} \right)^n C_n^{v+1}\left( \frac{c_1}{\sqrt{c_1^2-1}} \right)=
\frac{(-1)^n}{z_3^n}
C_n^{v+1}\left( \frac{y_1y_2}{\sqrt{y_1^2+y_2^2-1}} \right)
\left( y_1^2+y_2^2-1 \right)^{\frac{n}{2}}
\end{equation*}
Используя разложение для полиномов Гегенбауэра ( при $ \sigma=0 $ или $ \sigma=1 $  ) \\ \cite[Гл.4, п.4.6.2(6)]{Prudnikov_2}  
\begin{multline*}
C_{2n+\sigma}^{v+1}\left( \frac{y_1y_2}{\sqrt{y_1^2+y_2^2-1}} \right)
\left( y_1^2+y_2^2-1 \right)^{\frac{2n+\sigma}{2}}=
\frac{\Gamma( v+1 )\Gamma( 2n+2v+2+\sigma )}{2^{2n+\sigma}}
\times
\\
\times
\sum_{k=0}^{n}\frac{(2k+\sigma)!\,(2k+\sigma+v+1)}{(n-k)!\,\Gamma( 2k+2v+2+\sigma )\Gamma( n+k+v+\sigma+2 )}
C_{2k+\sigma}^{v+1}(y_1)C_{2k+\sigma}^{v+1}(y_2)
\end{multline*}
находим из \eqref{1_f3}, после объединения сумм по полиномам Гегенбауэра
\begin{multline*}
\frac{e^{-i\pi\mu}Q_{v}^{\mu}(c_1+c_2)}{\sqrt{(c_1+c_2)^2-1}^{\mu}}=
\frac{2^{\mu+2v+1}\Gamma( v+1 )^2}{(2c_2)^{v+\mu+1}}
\times
\\
\times
\sum_{n=0}^{\infty}
\frac{(-1)^n n!\,(n+v+1)}{\Gamma( n+2v+2 )}
C_n^{v+1}(y_1)C_n^{v+1}(y_2)
\sum_{k=0}^{\infty}\frac{\Gamma( v+\mu+1+2k+n )}{\Gamma( v+2+n+k )k!}\frac{1}{(2z_3)^{2k+n}}
\end{multline*} 
где $ c_1,c_2,y_1,y_2 $ выражаются из \eqref{1_f4} и \eqref{1_f5}.
Последнюю сумму по $ k $ представим в виде функции Лежандра второго рода из \eqref{1_f1}, а полиномы Гегенбауэра в $ z_k $ -переменных, и в итоге окончательно получим \eqref{1_f6}

Из \cite[Гл.3, п.3.6]{Beitman_1} функцию Лежандра второго рода $ Q_v^{n+v+1}(z) $ при целых $ n $ можно представить в виде конечных сумм. В частности
\begin{equation*}
Q_v^{n+v+1}(iz)=(-1)^n e^{i\tfrac{\pi}{2}(v+1)}\,\frac{2^v n!\,\Gamma( v+1 )}{\sqrt{z^2+1}^{v+1}}C_n^{v+1}\left( \frac{z}{\sqrt{z^2+1}} \right), \quad
(z>0,\quad n=0,1,2,\ldots)
\end{equation*}
Сопоставляя это выражение с \eqref{1_f6} можно получить теорему сложения следующего вида
\begin{multline} \label{1_f7}
\frac{e^{-i\,\pi\,\mu}Q_{v}^{\mu}(z)}{\sqrt{z^2-1}^{\mu}}=
\sqrt{\frac{8}{\pi}}\frac{e^{-i\pi(v+1)}}{\left(\sqrt{z_1^2+1}\sqrt{z_2^2+1} \right)^{\mu}}
\times
\\
\times
\sum_{n=0}^{\infty}
\frac{(-1)^n \, (n+v+1)}{\Gamma(n+2v+2)\,n!}
Q_v^{n+v+1}(iz_1)Q_v^{n+v+1}(iz_2)
\frac{e^{-i\,\pi\,\left( \mu -\frac{1}{2} \right)}Q_{n+v+\frac{1}{2}}^{\mu-\frac{1}{2}}(z_3)}{\sqrt{z_3^2-1}^{\mu -\frac{1}{2}}}
\end{multline}
\[ z=z_1z_2+z_3\sqrt{z_1^2+1}\sqrt{z_2^2+1},\quad (z_1>0,\,z_2>0,\,z_3>1)
\]
В частности при $ \mu=0 $ из \cite[Гл.3, п.3.6]{Beitman_1} можно также получить
\begin{equation*}
\lim_{\mu=0}
\left. 
\frac{e^{-i\,\pi\,\left( \mu -\frac{1}{2} \right)}Q_{n+v+\frac{1}{2}}^{\mu-\frac{1}{2}}(z_3)}{\sqrt{z_3^2-1}^{\mu -\frac{1}{2}}}
\right|_{z_3=\cosh\alpha}=-\sqrt{\frac{\pi}{2}}\,\frac{e^{-\alpha(n+v+1)}}{n+v+1}, \quad (\alpha>0)
\end{equation*}
в котором \eqref{1_f7} при целых $ v=0,1,2,\ldots $ будет соответствовать теореме сложения для функции Лежандра второго рода \cite[Гл.8, 8.795(3)]{Gradshteyn_Ryzhik}

\end{document}